\documentclass[11pt,a4paper]{amsart}
\usepackage[T1]{fontenc}
\usepackage{amsmath,amssymb,amsthm}
\usepackage{enumerate}
\newcommand\esp[2]{{\mathbb E}_{#1}\left[{#2}\right]}
\def\E{{\mathbb E}}
\def\det{{\text{det}}}

\def\E{{\mathbb E}}
\def\car{{\mathbf 1}}
\def\AA{{\mathfrak A}}
\def\P{{\mathbf P}}
\def\R{{\mathbf R}}
\def\N{{\mathbf N}}
\def\FC{{\mathcal{F}C^{\infty}_b}}
\def\BB{{\mathfrak B}}
\newcommand{\Det}{\operatorname{Det}}
\renewcommand{\d}{\text{ d}}
\newcommand{\Tr}{\operatorname{trace}}
\newcommand{\trace}{\operatorname{trace}}
\newcommand{\Diff}{\operatorname{Diff}}
\newcommand{\Jac}{\operatorname{Jac}}
\newcommand{\adj}{\operatorname{Adj}}
\newcommand{\I}{\operatorname{I}}
\newcommand{\divergence}{\operatorname{div}}
\newcommand{\per}{\operatorname{per}}
\newtheorem{Theorem}{Theorem}
\newtheorem{Lemma}{Lemma}
\newtheorem{Definition}{Definition}
\newtheorem{Remark}{Remark}

\newtheorem{Corollary}{Corollary}
\newtheorem{hyp}{Hypothesis}
\usepackage[latin1]{inputenc}
\usepackage{amsfonts}
\author{I. Camilier}
\address{
 Institut TELECOM, TELECOM ParisTech, CNRS LTCI\\
 Paris, France } 
\email{Isabelle.Camilier@telecom-paristech.fr}
\author{L. Decreusefond}
\address{
 Institut TELECOM, TELECOM ParisTech, CNRS LTCI\\
 Paris, France } 
\email{Laurent.Decreusefond@telecom-paristech.fr}

\title[Quasi-invariance of determinantal processes]{Quasi-invariance
  and integration by parts  for determinantal and permanental processes}
\date{\today}
\begin{document}
\def\K{{\text{\textbf{K}}}} \def\J{{\text{\textbf{J}}}}

\keywords{Determinantal point processes, Malliavin calculus,
  permanental point processes, point processes, integration by parts
  formula, random thinning.}  

\begin{abstract}
  Determinantal and permanental processes are point processes with a
  correlation function given by a determinant or a permanent. Their
  atoms exhibit mutual attraction of repulsion, thus these processes
  are very far from the uncorrelated situation encountered in Poisson
  models. We establish a quasi-invariance result : we show that if
  atoms locations are perturbed along a vector field, the resulting
  process is still a determinantal (respectively permanental) process,
  the law of which is absolutely continuous with respect to the
  original distribution. Based on this formula, following Bismut
  approach of Malliavin calculus, we then give an integration by parts
  formula.
\end{abstract}
\maketitle

\section{Motivations}
Point processes are widely used to model various phenomena, such as
arrival times, arrangement of points in space, etc. It is thus
necessary to know into details as large a catalog of point processes
as possible. The Poisson process is one example which has been widely
studied for a long time. Our motivation is to study point processes
that generate a more complex correlation structure, such as a
repulsion or attraction between points, but still remain simple enough
so that their mathematical properties are analytically tractable.
Determinantal and permanental point processes hopefully belong to this
category. They were introduced in~\cite{MR0380979} in order to
represent configurations of fermions and bosons. Elementary particles
belong exclusively to one of these two classes.  Fermions are
particles like electrons or quarks; they obey the Pauli exclusion
principle and hence the Fermi-Dirac statistics. The other sort of
particles are particles like photons which obey the Bose-Einstein
statistics. The interested reader can find in~\cite{MR2207648} an
illuminating account of the determinantal (respectively permanental)
structure of fermions (respectively bosons) ensemble.  A mathematical
unified presentation of determinantal/permanental point processes
(DPPP for short) was for the first time, introduced
in~\cite{MR2018415}. Let $\chi$ be the space of locally finite, simple
configurations on a Polish space $E$ and $K$ a locally trace-class
operator in $L^2(E)$ with a Radon measure $\lambda$. For any positive,
compactly supported $f$ and $\xi=\sum_j \delta_{x_j}\in \chi,$ the
$\alpha$-DPPP is the measure, $\mu_{\alpha,\, K, \, \lambda}$, on
$\chi$ such that
\begin{equation}\label{eq:1}
  \int_{\chi} e^{-\int f \d\xi }
  \d\mu_{\alpha,K,\lambda}(\xi) =\Det \left( \I+\alpha \sqrt{1-e^{-f}} \, K \, \sqrt{1-e^{-f}}\right)^{-\frac{1}{\alpha}},
\end{equation}
for the parameters $\alpha\in \AA=\{2/m;\, m\in \N\}\cup\{-1/m,\, ,\in
\N\},$ where $\N$ is the set of positive integers. The values
$\alpha=-1$ and $\alpha=1$ correspond to determinantal and permanental
point processes respectively. Starting from~\eqref{eq:1}, existence of
$\alpha$-DPPP is still a challenge as explained in
\cite{math.PR/0002099}. Actually, existence is (not easily) proved for
$\alpha=\pm 1$ and DPPP for other values of $\alpha$ are constructed
as superposition of these basic processes.  DPPP recently regained
interest because they have strong links with the spectral theory of
random matrices~\cite{math-ph0510038,math.PR/0002099}: for instance,
eigenvalues of matrices in the Ginibre ensemble a.s. form a
determinantal configuration.  DPPP also appear in polynuclear growth
\cite{MR2018275,math.PR/0304368}, non intersecting random walks,
spanning trees, zero set of Gaussian analytic functions (see
\cite{MR2216966} and references there in), etc. Mathematically
speaking, a few of their properties are known. The most complete
references to date are, to the best of our knowledge, \cite{MR2216966}
and \cite{MR2018415}. Gibbsianness of DPPP, i.e., local absolute
continuity of $\mu_{\alpha, K,\lambda}$ with respect to the
distribution of a Poisson process, was investigated in several papers
by Yoo \cite{MR2209150,MR2122549}. The conclusion of all these studies
seems to be that DPPP are rather hard to describe and analyze, their
properties being highly dependent of the kernel and its eigenvalues.

Our aim is to investigate further some of the stochastic properties of
$\alpha$-DPPP. In the spirit of \cite{MR2108363}, we are interested in
the differential calculus associated to these processes. In
\cite{MR2108363}, it is shown that a somewhat canonical Dirichlet form
associated to DPPP is closed. We here address the problem within the
point of view of Malliavin calculus. To date, Malliavin calculus for
point processes has been developed namely for Poisson processes
(\cite{bichteler83,bass86,MR2535466,MR99d:58179,MR2531026,Decreusefond:2006wr})
and some of their extensions: Gibbs processes \cite{MR99f:58219},
marked processes \cite{MR2002f:60098}, filtered Poisson processes
\cite{Decreusefond:2006wr}, cluster processes \cite{Bogachev:2007ww}
and Lévy processes \cite{MR2292579,Di-Nunno:2004no}. There exist three
approaches to construct a Malliavin calculus framework for point
processes: one based on white noise analysis, one based on a
difference operator and chaos decomposition and one which relies on
quasi-invariance of the law of Poisson process with respect to some
perturbations. This is the last track we follow here since neither the
white noise framework nor the chaos decomposition exist so far.

We first show that the action of a diffeomorphism of $E$ into itself
onto the atoms of a DPPP yields another DPPP, the law of which is
absolutely continuous with the distribution of the original process; a
property usually known as quasi-invariance. Then, following the lines
of proof of \cite{MR99d:58179,bismut83,Bogachev:2007ww}; we can derive
an integration by parts formula for the differential gradient as
usually constructed on configuration spaces. This gives another proof
of the closability of the Dirichlet form canonically associated to a
DPPP as in \cite{MR2108363}.

This paper is organized as follows. In part \ref{preliminaries}, we
give definitions concerning point processes and
$\alpha$-de\-ter\-mi\-nan\-tal point processes. In part
\ref{Quasi-invariance}, we prove the quasi-invariance for
$\alpha$-DPPP. Then, in Section \ref{Integration by parts}, we compute
the integration by parts formula. We begin by determinantal point
processes and then extend to $\alpha$-determinantal point
processes. Permanental processes are analyzed on the same basis.

\section{Preliminaries}\label{preliminaries}
\subsection{Point processes}
We remind here some properties of point processes we refer to
\cite{MR1950431,kallenberg83} for more details.  Let $E$ be a Polish
space and $\lambda$ a Radon measure on $(E,\, \mathcal{B})$, the Borel
$\sigma$-algebra on $E.$ By $\chi$ we denote the space of all locally
finite configurations on $E$:
\begin{equation*}\label{}
  \chi= \{ \xi\subset E: |\xi \cap \Lambda|< \infty
  \text{ for any compact } \Lambda \subset E  \},
\end{equation*}
where $|A|$ is the cardinality of a set $A$. Hereafter we identify a
locally finite configuration $\xi,$ defined as a set, and the atomic
measure $\sum_{x \in \xi}\delta_x$. The space $\chi$ is then endowed
with the vague topology of measures and $\mathcal{B}(\chi)$ denotes
the corresponding Borel $\sigma$-algebra. For any measurable
nonnegative function $f$ on $E$, we denote equivalently:
\begin{equation*}
  \left\langle f,\xi \right\rangle =
  \sum_{x\in \xi}f(x)=\int{f \d\xi}.
\end{equation*}
We also denote by $\chi_{0}=\left\lbrace \alpha \in \chi, \mid
  \alpha(E) \mid < \infty \right\rbrace$ the set of all finite
configurations in $\chi$ and $\chi_{0}$ is equipped with the
$\sigma$-algebra $\mathcal{B}(\chi_0)$.
The restriction of a configuration $\xi$ to a compact $\Lambda\subset
E$, is denoted by $\xi_\Lambda$. We introduce the set
$\chi_{\Lambda}=\left\lbrace \xi \in \chi, \xi(E\backslash\Lambda) =0
\right\rbrace$. Then for any integer $n$, we denote by
$\chi_{\Lambda}^{(n)}=\left\lbrace \xi \in \chi, \xi(\Lambda) =n
\right\rbrace$, the set of all configurations in with $n$ points in
$\Lambda$.  Note that we have
$\chi_{\Lambda}=\bigcup_{n=0}^{\infty}\chi_{\Lambda}^{(n)}$.
\begin{Definition}
  A random point process is a triplet $(\chi,\mathcal{B}(\chi),\mu)$,
  where $\mu$ is a probability measure on $(\chi,\mathcal{B}(\chi))$.
\end{Definition}
Every measure $\mu$ on the configuration space $\chi$ can be
characterized by its Laplace function, that is to say for any
measurable non-negative function $f$ on $E$:
\begin{equation*}
  f\longmapsto \mathbb{E}_{\mu}[e^{-\int{f \d\xi}}]=\int_{\chi}e^{-\int{f \d\xi}}\d\mu(\xi).
\end{equation*}
For instance, let $\pi_{\sigma}$ denote the Poisson measure on $(\chi,
\mathcal{B}(\chi))$ with intensity measure $\sigma$. Then its Laplace
transform is, for any measurable non-negative function $f$:
\begin{equation*}
  \int_{\chi}e^{-\int{f \d\xi}}\d\pi_{\sigma}(\xi)=\exp\left(\int_E(1-e^{-f(x)})\d\sigma(x)\right)
\end{equation*}
Another way to describe the distribution of a point process is to give
the probabilities $\P(|\xi_{\Lambda_k}|=n_k, \, 1 \le k \le n )$ for
any $n$ and any mutually disjoints Borel subsets of $\Lambda$,
$\Lambda_1,\cdots,\, \Lambda_k$, $1 \le k \le n$.  For instance, the
Poisson measure $\pi_{\sigma}$ with intensity measure $\sigma$ can be
defined in this way as:
\begin{equation*}
  \P(|\xi_{\Lambda_k}|=n_k, \, 1 \le k \le n )=\prod_{k=1}^{n}e^{-\sigma(\Lambda_k)}\frac{\sigma(\Lambda_k)^{n_k}}{n_k!}.
\end{equation*}
But in many cases, specifying the joint distribution of the $\xi(D)$'s
is not that simple. It is then easier to describe the distribution of
a point process by its correlation functions.
\begin{Definition}
  A locally integrable function $\rho_n: E^n\rightarrow\R_{+}$ is the
  $n$-point correlation function of $\mu$ if for any disjoint bounded
  Borel subsets $\Lambda_1,\, \cdots,\, \Lambda_m$ of $E$ and $n_i \in
  \N$, $\sum_{i=1}^{m}n_i=n$:
  \begin{equation*}
    \mathbb{E}_{\mu}\left[\prod_{i=1}^{m}{\dfrac{|\xi_{\Lambda_i}|!}{(|\xi_{\Lambda_i}|-n_i)!}}\right]=\int_{\Lambda_1^{n_1}\times\ldots\times\Lambda_m^{n_m}}{\rho_n(x_1,\,
      \cdots,\, x_n)\d\lambda(x_1)\ldots\d\lambda(x_n)},
  \end{equation*}
  where $\mathbb{E}_{\mu}$ denotes the expectation relatively to
  $\mu.$
\end{Definition}
For example, if $m=1$ and $n_1=n$, the formula becomes:
\begin{multline*}
  \esp{\mu}{\frac{|\xi_\Lambda|!}{(|\xi_{\Lambda}|-n)!}}
  =\esp{\mu}{|\xi_{\Lambda}|\, (|\xi_\Lambda|-1)\ldots(|\xi_\Lambda|-n+1)}\\
  =\int_{\Lambda^n}{\rho_n(x_1,\, \cdots,\,
    x_n)\d\lambda(x_1)\ldots\d\lambda(x_n)}.
\end{multline*}
We recognize here the $n$-th factorial moment of $|\xi_\Lambda|$.  In
particular:
\begin{equation*} \esp{\mu}{|\xi_\Lambda|}
  =\int_{\Lambda}{\rho_1(x)\d\lambda(x)},
\end{equation*}
i.e., $\rho_1$ is the mean density of particles. More generally, the
function $\rho_n$ has the following interpretation: $\rho_n(x_1,\,
\cdots,\, x_n)\d\lambda(x_1)\ldots\d\lambda(x_n)$ is approximately the
probability to find a particle in each one of the $[x_i,\, x_i+\d
\lambda(x_i)]$, $i=1,\cdots, n$.  A third way to define a point
process proceeds via the Janossy densities. Denote by
$\pi_{n,\Lambda}(x_1,\, \cdots,\, x_n)$ the density (assumed to exist)
with respect to $\lambda^{\otimes n}$ of the joint distribution of
$(x_1,\, \cdots,\, x_n) $ given that there are $n$ points in
$\Lambda$.
\begin{Definition}
  The density distributions or Janossy densities of a random process
  $\mu$ are the measurable functions $j^n_{\Lambda}$ such that:
  \begin{eqnarray*}
    j^n_{\Lambda}(x_1,\, \cdots,\, x_n) &=& n!\, \mu(\xi(\Lambda)=n)\, \pi_{n,\Lambda}(x_1,\, \cdots,\, x_n)
    \text{ for }n \in \N,\\
    j^0_{\Lambda}(\emptyset) &=& \mu(\xi(\Lambda)=0).
  \end{eqnarray*}
\end{Definition}
Hence, the Janossy density $j^n_{\Lambda}(x_{1},\, \cdots,\, x_{n})$
appears as the probability density that there are exactly $n$ points
in $\Lambda$ located around $%
x_{1},\, \cdots,\, x_{n}$, and no points anywhere else. For
$n=0$, $j^0_{\Lambda}(\emptyset)$ is the probability that there is no
point in $\Lambda$. For $n \geq 1$, the Janossy densities satisfy the
following properties:
\begin{itemize}
\item Symmetry: \begin{equation*} j^n_{\Lambda}\left( x_{\sigma(1)},\,
      \cdots,\, x_{\sigma(n)}\right) =j_{n,\Lambda}\left( x_{1},\,
      \cdots,\, x_{n}\right),
  \end{equation*}
  for every permutation $\sigma$ of $\left\lbrace 1,\, \cdots,n
  \right\rbrace$.
\item Normalization constraint. For each compact $\Lambda$:
  \begin{equation*}
    \sum_{n=0}^{+\infty}\int_{\Lambda^n}{
      \frac{1}{n!}\, j^n_{\Lambda}\left( x_{1},\, \cdots,\, x_{n}\right)
      \d\lambda(x_1)\ldots\d\lambda(x_n) }=1.
  \end{equation*}
\end{itemize}
It is clear that the $\rho_n$'s, $j_n$'s, $\mu$ should satisfy some
relationships. We will not dwell on that here (see the references
cited above), we just mention the relation between $\mu$ and
$j^n_{\Lambda},$ which is:
\begin{multline}\label{janossy}
  \int_{\chi}{f\left( \xi \right) \d\mu(\xi) }=\\
  \sum_{n=0}^{+\infty}\frac{1}{n!}\int_{\Lambda^n}{f(x_1,\, \cdots,\,
    x_n)\ j^n_{\Lambda}\left( x_{1},\, \cdots,\, x_{n}\right)
    \d\lambda(x_1)\ldots\d\lambda(x_n) }.
\end{multline}

\subsection{Fredholm determinants}
For details on this part, we refer to \cite{0635.47002,MR2154153}. For
any compact $\Lambda \subset E$, we denote by $L^2(\Lambda,\,
\lambda)$ the set of functions square integrable with respect to the
restriction of the measure $\lambda$ to the set $\Lambda$. This
becomes a Hilbert space when equipped with the usual norm:
\begin{equation*}
  \|f\|^2_{L^2(\lambda,\Lambda)}=\int_{\Lambda}|f(x)|^2 \d\lambda(x).
\end{equation*}
For $\Lambda$ a compact subset of $E$, $P_{\Lambda}$ is the projection
from $L^2(E)$ onto $L^2(\Lambda)$, i.e., $P_\Lambda f=f\car_\Lambda.$
The operators we will deal with are special cases of the general
category of continuous maps from $L^2(E,\, \lambda)$ into itself.
\begin{Definition}
  A map $T$ from $L^2(E)$ into itself is said to be an integral
  operator whenever there exists a measurable function, we still
  denote by $T$, such that
  \begin{equation*}
    Tf(x)=\int_E T(x,\, y) f(y)\d \lambda(y).
  \end{equation*}
  The function $T$ is called the kernel of $T$.
\end{Definition}
\begin{Definition}
  Let $T$ be a bounded map from $L^2(E,\, \lambda)$ into itself. The
  map $T$ is said to be trace-class whenever for one complete
  orthonormal basis (CONB for short) $(h_n,\, n\ge 1)$ of $L^2(E,\,
  \lambda)$,
  \begin{equation*}
    \sum_{n\ge 1} |(Th_n,\, h_n)_{L^2}|\text{ is finite.}
  \end{equation*}
  Then, the trace of $T$ is defined by
  \begin{equation*}
    \trace(T)= \sum_{n\ge 1} (Th_n,\, h_n)_{L^2}.
  \end{equation*}
\end{Definition}
It is easily shown that the notion of trace does not depend on the
choice of the CONB. Note that if $T$ is trace-class then $T^n$ also is
trace-class for any $n\ge 2.$
\begin{Definition}
  \label{def:fredholm_determinant}
  Let $T$ be a trace class operator. The Fredholm determinant of
  $(\I+T)$ is defined by:
  \begin{equation*}
    \Det(\I+T)=\exp\left(\sum_{n=1}^{+\infty}\frac{(-1)^{n-1}}{n}\Tr (T^n)\right),
  \end{equation*}
  where $\I$ stands for the identity operator.
\end{Definition}
The practical computations of fractional power of Fredholm
determinants involve the so-called $\alpha$-determinants, which we
introduce now.
\begin{Definition}
  For a square matrix $A=\left( a_{ij}\right)_{i,j=1\ldots n} $ of
  size $n \times n$, the $\alpha$-determinant $\det_{\alpha}\, A$ is
  defined by:
  \begin{equation*}
    \det_{\alpha}\, A=\sum_{\sigma\in\Sigma_n}\alpha^{n-\nu\left(
        \sigma\right) }\prod_{i=1}^{n}a_{i\sigma\left( i\right) },
  \end{equation*}
  where the summation is taken over the symmetric group $\Sigma_n$,
  the set of all permutations of $\left\lbrace 1, 2,\, \cdots, n
  \right\rbrace$ and $\nu(\sigma)$ is the number of cycles in the
  permutation ~$\sigma$.
\end{Definition}
This is actually a generalization of the well-known determinant of a
matrix. Indeed, when $\alpha=-1$, $\det_{-1}\, A$ is the usual
determinant $\det\, A.$ When $\alpha=1$, $\det_{1} \, A$ is the
so-called permanent of $A$ and for $\alpha=0$, $\det_{0}\, A=\prod_{i}
a_{ii}$.  We can then state the following useful theorem (see
\cite{MR2018415}):
\begin{Theorem}
  \label{thm:developpement_det_alpha}
  For a trace class integral operator $T$, if $\parallel \alpha T
  \parallel < 1 $, we have:
  \begin{equation*}
    \Det (\I-\alpha T)^{-\frac{1}{\alpha}}=
    \sum_{n=0}^{+\infty}\frac{1}{n!}\int_{\Lambda^n}{\det_{\alpha}\, (T(x_i,\, x_j))_{1\le i,j\le n}}\d\lambda(x_1)\ldots\d\lambda(x_n).
  \end{equation*}
  If $\alpha \in \left\lbrace -1/m; m \in \N \right\rbrace $, this is
  true without the condition $\parallel \alpha T \parallel < 1 $.
\end{Theorem}
\subsection{Determinantal-permanental point processes}
\label{sec:alpha-dppp}
The following set of hypothesis is of constant use.
\begin{hyp}\label{hyp:condition_K} The Polish space $E$ is equipped
  with a Radon measure $\lambda$. The map $K$ is an Hilbert-Schmidt
  operator from $L^2(E,\, \lambda)$ into $L^2(E,\, \lambda)$ which
  satisfies the following conditions:
  \begin{enumerate}[i)]
  \item \label{item:1} $K$ is a bounded symmetric integral operator on
    $L^2(E,\, \lambda)$, with kernel $K(.,.)$, i.e., for any $x\in E$,
    \begin{equation*}
      K f(x)=\int_{E}{K(x,y)f(y)\d\lambda(y)}.
    \end{equation*}
  \item \label{item:2} The spectrum of $K$ is included $[ 0,\, 1[$.
  \item \label{item:3} The map $K$ is locally of trace class, i.e.,
    for all compact $\Lambda\subset E$, the restriction
    $K_{\Lambda}=P_{\Lambda}K P_{\Lambda}$ of $K$ to $L^2(\Lambda)$ is
    of trace class. 
  \end{enumerate}
\end{hyp}
For a real $\alpha \in [-1,1]$ and a compact subset $\Lambda\subset
E$, the map $J_{\Lambda,\alpha}$ is defined by:
\begin{equation*}
  J_{\Lambda,\alpha} =\left(\I+\alpha
    K_{\Lambda}\right)^{-1} K_{\Lambda},
\end{equation*}
so that we have:
\begin{equation*}
  \left(\I+\alpha K_{\Lambda} \right) \left(\I- \alpha
    J_{\Lambda,\alpha}\right) =\I.
\end{equation*}
For any compact $\Lambda$, the operator $J_{\Lambda,\alpha}$ is also a
trace class operator in $L^2(\Lambda,\, \lambda)$.  In the following
theorem, we define $\alpha$-DPPP with the three equivalent
characterizations: in terms of their Laplace transforms, Janossy
densities and correlation functions. The theorem is also a theorem of
existence, a problem which as said above is far from being trivial.
\begin{Theorem}[See \protect{\cite{MR2018415}}]\label{thm:existence}
  Assume Hypothesis~\ref{hyp:condition_K} is satisfied.  Let $\alpha
  \in \AA$. There exists a unique probability measure
  $\mu_{\alpha,\,K,\,\lambda}$ on the configuration space $\chi$ such
  that, for any nonnegative bounded measurable function $f$ on $E$
  with compact support, we have:
  \begin{align}
    \E_{\mu_{\alpha,\,K,\,\lambda}}\left[ e^{-\int f \d\xi }\right]&=
    \int_{\chi} e^{-\int f \d\xi }
    \d\mu_{\alpha,\,K,\,\lambda}(\xi)\notag\\
    & =\Det \left( \I+\alpha K
      [1-e^{-f}]\right)^{-\frac{1}{\alpha}},\label{eq:3}
  \end{align}
  where $K[1-e^{-f}]$ is the bounded operator on $L^2(E)$ with kernel
  :
  \begin{equation*}
    (K[1-e^{-f}])(x,y)=  \sqrt{1-\exp(-f(x))} K(x,y)\, \sqrt{1-\exp(-f(y))}.
  \end{equation*}
  This means that for any integer $n$ and any $(x_1,\cdots,\, x_n) \in
  E^n,$ the correlation functions of $\mu_{\alpha,\, K,\, \lambda}$
  are given by:
  \begin{equation*}
    \rho_{n,\,\alpha,\, K}(x_1,\cdots,\,\, x_n)=\det_{\alpha}\left( K\left( x_i,\, x_j\right)
    \right)_{1\le i,j\le n},
  \end{equation*}
  and for $n=0$, $\rho_{0,\,\alpha,\,K}(\emptyset)=1$.  For any
  compact $\Lambda\subset E,$ the operator $J_{\Lambda,\, \alpha}$ is
  an Hilbert-Schmidt, trace class operator, whose spectrum is included
  in $[0, +\infty[$.  For any $n \in \N$, any compact $\Lambda\subset
  E$, and any $(x_1,\cdots,\, x_n) \in \Lambda^n$ the $n$-th Janossy
  density is given by:
  \begin{equation}\label{eq:4}
    j_{\Lambda,\alpha,K_{\Lambda}}^{n}\left( x_1,\, \cdots,\, x_n\right)
    =\Det\left( \I+\alpha
      K_{\Lambda}\right)^{-1/\alpha}\det_{\alpha}\left(J_{\Lambda,\,\alpha} (x_i,\, x_j)\right) _{1\le i,j \le n}.
  \end{equation}
  For $n=0$, we have
  \begin{math}
    j_{\Lambda,\alpha,K_{\Lambda}}^{n}\left( \emptyset\right)
    =\Det\left( I+\alpha K_{\Lambda}\right)^{-1/\alpha}.
  \end{math}
\end{Theorem}
For $\alpha=-1$, such a process is called a determinantal process
since we have, for any $n \geq 1$:
\begin{equation*}
  \rho_{n,-1,K}(x_1,\cdots,\,\, x_n)=\det(K(x_i,\, x_j))_{1 \le i,j \le n}.
\end{equation*}
For $\alpha=1$, such a process is called a permanental process, since
we have, for any $n \geq 1$:
\begin{equation*}
  \rho_{n,1,K}(x_1,\cdots,\, x_n)=\sum_{\pi\epsilon\Sigma}\prod_{i=1}^{n}K(x_{i},\, x_{\pi(i)})=\text{per}\, (K(x_i,\, x_j))_{1 \le i,j \le n}.
\end{equation*}
For any bounded function $g:E\rightarrow \R^+$, and any integral
operator $T$ of kernel $T(x,y)$, we denote by $T[g]$ the integral
operator of kernel:
\begin{equation*}
  T[g](x,y)\rightarrow \sqrt{g(x)} T(x,y) \sqrt{g(y)}.
\end{equation*}
For calculations, it will be convenient to use the following lemma:
\begin{Lemma}[see \cite{MR2018415}]
  \label{lem:calcul}
  Let $\Lambda$ be a compact subset of $E$ and $f:E\rightarrow [0,
  +\infty),$ measurable with supp$(f) \in \Lambda$:
  \begin{align*}
    \Det\left( \I+ \alpha K_{\Lambda}[1-e^{-f}]\right)
    ^{-1/\alpha}
    &=\Det\left( \I+\alpha K_{\Lambda}\right) ^{-1/\alpha}\Det\left(
      \I- \alpha J_{\Lambda,\alpha} [e^{-f}] \right)^{-1/\alpha}.
  \end{align*}
\end{Lemma}
By differentiation into the Laplace transform, it is possible to
compute moments of $\int f \d\xi$ for any deterministic $f$. We
obtain, at the first order:
\begin{Theorem}[see
  \cite{MR2018415}] \label{thm:esperance_et_variance} For any
  non-negative function $f$ defined on $E$, we have
  \begin{equation*}
    \esp{}{\int_\Lambda f\d\xi}=\int_\Lambda f(x)K(x,\, x)\d\lambda(x)=\trace(K_\Lambda[f]).
  \end{equation*}
\end{Theorem}
It is worth mentioning how the existence of $\alpha$-DPPP is
established. For $\alpha=-1,$ there is a non trivial work (see
\cite{MR2018415,math.PR/0002099} and references therein) to show that
the Janossy densities satisfy the positivity condition so that a point
process with these densities does exist. For $\alpha=-1/m$, it is
sufficient to remark from (\ref{eq:3}) that the superposition of $m$
independent determinantal point processes of kernel $K/m$ is an
$\alpha$-DPPP fo kernel $K$. The point is that $K/m$ satisfies Hypothesis
\ref{hyp:condition_K}, in particular that its spectrum is strictly
bounded by $1/m< 1,$ since $m> 1$.  For $\alpha=2$, a
$2$-permanental point process is in fact a Cox process based on a
Gaussian random field. We know for sure that there exists $X$ a
centered Gaussian random field on $E$ such that:
\begin{equation}
  \label{eq:5b}
  \esp{\P}{\int_\Lambda X^2(x)\d \lambda(x)}=\trace(K_\Lambda),
\end{equation}
for any compact $\Lambda\subset E$ and
\begin{equation}
  \label{eq:6}  \esp{\P}{X(x)X(y)}=K(x,\, y) \ \lambda\otimes \lambda
  \text{ a.s.}, 
\end{equation}
where $\P$ is the probability measure on the probability space
supporting $X$.  Then the Cox process of random intensity $X^2(x)\d
\lambda(x)$ has the same distribution as $\mu_{2,K,\lambda}.$ Indeed,
it follows from the formula:
\begin{equation*}
  \esp{\P}{\exp\left(-\int (1-e^{-f(x)}) X^2(x)\d \lambda(x)\right)}=\Det(I+2(1-e^{-f})K)^{-1/2}.
\end{equation*}
Thus, any $2/m$-permanental point process is the superposition of $m$
independent $2$-permanental point processes with kernel $K/m$.

Poisson process can be obtained formally as extreme case of
$1$-permanental process with a kernel $K$ given by $K(x,\,
y)=\car_{\{x=y\}}$. Of course, this kernel is likely to be null almost
surely with respect to $\lambda\otimes \lambda$; nonetheless, it
remains that replacing formally this expression in (\ref{eq:3}) yields
the Laplace transform of a Poisson process of intensity
$\lambda$. Another way to retrieve a Poisson process is to let
$\alpha$ go to $0$ in (\ref{eq:3}). With the above constructions, this
means that a Poisson process can be viewed as an infinite
superposition of determinantal or permanental point processes.
\begin{Theorem}
  When $\alpha$ tends to $0$, $\mu_{\alpha,K,\lambda}$ converges
  narrowly to a Poisson measure of intensity $K(x,\, x)\d\lambda(x).$
\end{Theorem}
\begin{proof}
  For any nonnegative $f$, for any $n\ge 1$,
  \begin{equation*}
    0\le
    \trace\left((K_{\Lambda}[1-e^{-f}])^n\right)\le
    \trace\left(K_{\Lambda}[1-e^{-f}]\right),
  \end{equation*}
  hence,
  \begin{multline}\label{eq:5}
    \int_{\chi}\exp\left(-\int f \ \d\xi
    \right)\d\mu_{\alpha,K_{\Lambda},\lambda}(\xi)
    {=\Det\left(\I+\alpha K_{\Lambda}[1-e^{-f}]\right)^{-1/{\alpha}}}\\
    {=\exp\left(-\frac{1}{\alpha}\sum_{n=1}^{+\infty}\frac{(-1)^{n-1}}{n}\alpha^n \Tr((K_{\Lambda}[1-e^{-f}])^n)\right)}\\
    {\xrightarrow{\alpha\rightarrow
        0}\exp\left(-\Tr(K_{\Lambda}(1-e^{-f}))\right)}=\int_E(1-e^{-f(x)})K_{\Lambda}(x,\,
    x)d\lambda(x).
  \end{multline}
  Thus, when $\alpha$ goes to $0$, the measure
  $\mu_{\alpha,K_{\Lambda},\lambda}$ tends towards a measure that we
  call $\mu_{0,K_{\Lambda},\lambda}$. According to (\ref{eq:5}),
  $\mu_{0,K_{\Lambda},\lambda}$ is a Poisson process with intensity
  $K_{\Lambda}(x,\, x)d\lambda(x)$.
\end{proof}
\section{Quasi-invariance}
\label{Quasi-invariance}
In this part we show the quasi-invariance property for any
$\alpha$-DPPP.  Let $\Diff _{0}(E)$ be the set of all diffeomorphisms
from $E$ into itself with compact support, i.e., for any $\phi\in
\Diff_0(E),$ there exists a compact $\Lambda$ outside which $\phi$ is
the identity map.  For any $\xi \in \chi$, we still denote by $\phi$
the map:
\begin{align*}
  \phi\, :\, \chi & \longrightarrow \chi\\
  \sum_{x \in \xi}\delta_{x}&\longmapsto \sum_{x \in
    \xi}{\delta_{\phi(x)}} .
\end{align*}
For any reference measure $\lambda$ on $E$, $\lambda_{\phi}$ denotes
the image measure of $\lambda$ by $\phi$.  For $\phi \in \Diff _0(E)$
whose support is included in $\Lambda$, we introduce the
isometry~$\Phi$,
\begin{align*}
  \Phi \, :\, L^2(\lambda_{\phi}, \Lambda)& \longrightarrow L^2(\lambda, \Lambda)\\
  f &\longmapsto f\circ \phi.
\end{align*}
Its inverse, which exists since $\phi$ is a diffeomorphism, is
trivially defined by $f\circ \phi^{-1}$ and denoted by
$\Phi^{-1}$. Note that $\Phi$ and $\Phi^{-1}$ are isometries, i.e.,
\begin{equation*}
  \langle \Phi \psi_1 ,\ \Phi \psi_2\rangle_{L^2({\lambda,\Lambda})}= \langle \psi_1 ,\ \psi_2\rangle_{L^2({\lambda_{\phi},\Lambda})}, 
\end{equation*}
for any $\psi_1$ and $\psi_2$ belonging to $L^2(\lambda, \Lambda)$.
We also set:
\begin{equation*}
  K^{\phi}_{\Lambda}=\Phi^{-1}K_{\Lambda}\Phi \text{ and }
  J^{\phi}_{\Lambda,\alpha}=\Phi^{-1}J_{\Lambda,\alpha}\Phi.
\end{equation*}
\begin{Lemma}
  \label{thm:operateurs}
  Let $\lambda$ be a Radon measure on $E$ and $K$ a map satisfying
  hypothesis \ref{hyp:condition_K}.  Let $\alpha \in \AA$. We have the
  following properties.
  \begin{enumerate}[a)]
  \item \label{item:4} $K^{\phi}_{\Lambda}$ and
    $J^{\phi}_{\Lambda,\alpha}$ are continuous operators from
    $L^2(\lambda_{\phi}, \Lambda)$ into $L^2(\lambda_{\phi},
    \Lambda)$.
  \item \label{item:5} $K^{\phi}_{\Lambda}$ is of trace class and
    $\Tr(K^{\phi}_{\Lambda})=\Tr(K_\Lambda)$.
  \item \label{item:8} $\Det(\I+\alpha
    K_{\Lambda}^{\phi})=\Det(\I+\alpha K_{\Lambda})$.
  \end{enumerate}
\end{Lemma}
\begin{proof}
  The first point is immediate according to the definition of an image
  measure.  Since $\Phi^{-1}$ is an isometry, for any $(\psi_n, \,
  n\in \N)$ a complete orthonormal basis of $L^2(\lambda, \Lambda)$,
  the family $(\Phi^{-1}\psi_n,\, ,\in \N)$ is a CONB of
  $L^2(\lambda_{\phi}, \Lambda)$. Moreover,
  \begin{align*}
    \sum_{n\ge 1} \left|\langle K^{\phi}_{\Lambda} \Phi^{-1}\psi_n,\,
      \Phi^{-1}\psi_n
      \rangle_{L^2(\lambda_{\phi}, \Lambda)}\right|&= \sum_{n\ge 1} \left|\langle \Phi^{-1} K \Phi \Phi^{-1}\psi_n,\, \Phi^{-1}\psi_n\rangle_{L^2(\lambda_{\phi}, \Lambda)}\right|\\
    &=\sum_{n\ge 1} \left|\langle \Phi^{-1} K\psi_n,\, \Phi^{-1}\psi_n\rangle_{L^2(\lambda_{\phi}, \Lambda)}\right|\\
    &=\sum_{n\ge 1} \left|\langle K \psi_n,\,
      \psi_n\rangle_{L^2(\lambda, \Lambda)}\right|.
  \end{align*}
  Hence, $K^{\phi}_\Lambda$ is of trace class and
  $\Tr(K^{\phi}_\Lambda)=\Tr(K_\Lambda)$.  Along the same lines, we
  prove that $ \Tr((K^{\phi}_\Lambda)^n)=\Tr(K_\Lambda^n)$ for any
  $n\ge 2.$ According to Definition \ref{def:fredholm_determinant},
  the Fredholm determinant of $K^\phi_{\Lambda}$ is well defined and
  point \ref{item:8}) follows.
\end{proof}
\begin{Theorem}
  \label{thm:kernel_form}
  Assume that $K$ is a kernel operator.  Then $K^{\phi}_{\Lambda},$ as
  a map from $L^2(\lambda_{\phi},\Lambda)$ into itself is a kernel
  operator whose kernel is given by $((x,\, y)\mapsto
  K_{\Lambda}(\phi^{-1}(x),\phi^{-1}(y))).$ An analog formula also
  holds for the operator $J_{\Lambda,\alpha}.$
\end{Theorem}
\begin{proof}
  On the one hand, for any function $f$, the operator
  $K^{\phi}_{\Lambda}$ from $L^2(\Lambda,\lambda_{\phi})$ into
  $L^2(\Lambda,\lambda_{\phi})$ is given by :
  \begin{equation*}
    K^{\phi}_{\Lambda}f(x)=\int_{\Lambda}K^{\phi}_{\Lambda}(x,z)f(z)d\lambda_{\phi}(z).
  \end{equation*}
  On the other hand, using the definition
  $K^{\phi}_{\Lambda}=\Phi^{-1}K_{\Lambda}\Phi$
  \begin{align*}
    K^{\phi}_{\Lambda}f(x)&=\Phi^{-1}K_{\Lambda}\Phi f(x)\\
    &=\int_{\Lambda}K_{\Lambda}(\phi^{-1}(x),y) f\circ \phi(y)d\lambda(y)\\
    &=\int_{\Lambda}K_{\Lambda}(\phi^{-1}(x),\phi^{-1}(z))
    f(z)d\lambda_{\phi}(z).
  \end{align*}
  The proof is thus complete.
\end{proof}

\begin{Lemma}
  \label{lem:chgt_mesure_operateur}
  Let $\rho\, :\, E\to \R$ be non negative and assume that
  $\d\lambda=\rho\d m$ for some other Radon measure on $E$. Let $K$
  satisfy Hypothesis \ref{hyp:condition_K}. Then, we have the
  following properties:
  \begin{enumerate}
  \item The map $K[\rho]$ is continuous from $L^2(m)$ into itself.
  \item The map $K[\rho]$ is locally trace class and $\trace(K
    _\Lambda[\rho])=\trace(K_\Lambda).$
  \item The measure $\mu_{\alpha, K, \lambda}$ is identical to the
    measure $\mu_{\alpha, K[\rho],m}.$
  \end{enumerate}
  That is to say, in some sense, we can ``transfer'' a part of the
  reference measure into the operator and vice-versa.
\end{Lemma}
\begin{proof}
  Remember that
  \begin{equation*}
    K[\rho] (x,y)=\sqrt{\rho(x)}\, K(x,y)\,\sqrt{\rho(y)}.
  \end{equation*}
  Hence
  \begin{equation*}
    K_{\Lambda}[\rho] f(x)=\sqrt{\rho(x)}\int_\Lambda
    K_{\Lambda}(x,y)\sqrt{\rho(y)}\d \lambda(y),
  \end{equation*}
  thus
  \begin{equation*}
    \int_\Lambda |K_{\Lambda}[\rho] f|^2 \d m=\int_\Lambda |K_{\Lambda} f|^2 \d \lambda,
  \end{equation*}
  and the first point follows.  Consider $(\psi_n,\, n \in \N)$, a
  CONB of $L^2(\lambda)$. Then $({\psi_n}{\sqrt{\rho}}, n \in \N)$ is
  a CONB of $L^2(m)$. Furthermore, we have:
  \begin{align*}
    \sum_{n\geq 1}\left|\langle K_{\Lambda}[\rho]
      \psi_n,\psi_n\rangle_{L^2(dm)}\right|&=\sum_{n\geq
      1}\left|\langle K_{\Lambda}\sqrt{\rho}
      \psi_n,\sqrt{\rho}\psi_n\rangle_{L^2(dm)}\right|\\
    &=\sum_{n\geq 1}\left|\left\langle K_{\Lambda}
        {\psi_n},{\psi_n}\right\rangle_{L^2(\lambda)}\right|.
  \end{align*}
  Therefore the operator $K_{\Lambda}[\rho]$ is of trace class and
  \begin{equation*}
    \Tr(K_{\Lambda}[\rho])=\Tr(K_{\Lambda}).
  \end{equation*}
  Similarly we can prove that for any $n \geq 2,$ we have
  $\Tr(K^n_{\Lambda}[\rho])=\Tr(K^n_{\Lambda})$. Then, using the
  definition of a Fredholm determinant, we have:
  \begin{equation*}
    \Det(\I+\alpha K_{\Lambda})=\Det(\I+\alpha K_{\Lambda}[\rho]).
  \end{equation*}
  The third point then follows from the characterization of
  $\mu_{\alpha, K[\rho],m}$ by its Laplace transform.
\end{proof}
The expression $\det_{\alpha}J_{\Lambda,\alpha}(x_i,\,\, x_j)_{1 \leq
  i,j \le n}$ is now denoted
$\det_{\alpha}J_{\Lambda,\alpha}(x_1,\cdots,\, x_n)$.  For any finite
random configuration $\xi=(x_1,\cdots,\, x_n)$, we call
$J_{\Lambda,\alpha}(\xi)$ the matrix with terms
$(J_{\Lambda,\alpha}(x_i,\, x_j), \ 1\le i,\, j\le n).$ First, remind
some results from \cite{MR99d:58179} concerning Poisson measures. For
any $\phi \in \Diff_0(E)$, we define $\phi^*\pi_{\lambda}$ as the
image of the Poisson measure $\pi_{\lambda}$ with intensity measure
$\lambda$ and $\lambda_{\phi}$ denotes the image measure of $\lambda$
by $\phi$.
\begin{Theorem}[\protect{\cite{MR99d:58179}}]
  \label{thm:invariance_poisson}
  For any $\phi \in \Diff _0(E)$, and a Poisson measure
  $\pi_{\lambda}$ with intensity $\lambda$:
  \begin{equation*} {\phi}^*\pi_{\lambda}=\pi_{\lambda_{\phi}}.
  \end{equation*}
  That is to say, for any $f$ nonnegative and compactly supported on
  $E$:
  \begin{equation}
    \mathbb{E}_{\pi_{\lambda}}\left[ e^{-\int{f\circ \phi
          \d\xi}}\right] =\exp\left(-\int 1-e^{-f} \d \lambda_\phi\right).
  \end{equation}
\end{Theorem}
We give the corresponding formula for $\alpha$-determinantal
measures. For any $\phi \in \Diff _0(E)$, we define
$\phi^*\mu_{\alpha,K_{\Lambda},\lambda}$ as the image of the measure
$\mu_{\alpha,K_{\Lambda},\lambda}$ under $\phi$. We prove below that
this image measure is an $\alpha$-DPPP the parameters of which are
explicitely known.
\begin{Theorem}
  \label{thm:changementdeloi}
  With the notations and hypothesis introduced above. For any $\phi
  \in \Diff _0(E)$, for any nonnegative function $f$ on $E$, for any
  compact $\Lambda\subset E$, we have:
  \begin{align}
    \mathbb{E}_{\mu_{\alpha,K_{\Lambda},\lambda}}\left[
      e^{-\int{f\circ \phi \d\xi}}\right] &
    =\mathbb{E}_{\mu_{\alpha,K_{\Lambda}^{\phi},\lambda_{\phi}}}\left[
      e^{-\int{f\d\xi}}\right]\label{eq:2}\\
    &=\Det(I+\alpha K_\Lambda^\phi[1-e^{-f}])^{-1/\alpha}.\notag
  \end{align}
  That is to say the image measure of $\mu_{\alpha,K,\lambda}$ by
  $\phi$ is an $\alpha$-determinantal process with operator $K^{\phi}$
  and reference measure $\lambda_\phi$.
\end{Theorem}
\begin{proof}
  According to Theorem \ref{thm:existence} and Lemma \ref{lem:calcul},
  we have for nonnegative $f$:
  \begin{multline*}
    \E_{\mu_{\alpha,K_{\Lambda},\lambda}} \left[ e^{- \int{f\circ \phi
          \ \d\xi}} \right]=
    \Det\left( \I+\alpha K_{\Lambda}[1-e^{-f\circ \phi} ] \right)^{-1/\alpha}\\
    =\Det\left( \I+\alpha K_{\Lambda}\right)^{-1/\alpha} \Det\left(
      \I-\alpha J_{\Lambda,\alpha} [e^{-f\circ
        \phi}]\right)^{-1/\alpha}.
  \end{multline*}
  According to Theorem \ref{thm:developpement_det_alpha}, we get
  \begin{multline*}
    \Det\left( \I-\alpha J_{\Lambda,\alpha} [e^{-f\circ \phi}]\right)^{-1/\alpha}\\
    \begin{aligned}
      &{=\sum_{n=0}^{+\infty}\frac{1}{n!}
        \int_{\Lambda^n}{\det_{\alpha}\, J_{\Lambda,\alpha}(x_1,\,
          \cdots,\, x_n)
          e^{-\sum_{i=1}^n f(\phi(x_i))}\d\lambda(x_1)\ldots \d\lambda(x_n)  }}\\
      &{=\sum_{n=0}^{+\infty}\frac{1}{n!}\int_{\Lambda^n}{\det_{\alpha}\,
          J^{\phi}_{\Lambda,\alpha}(x_1,\, \cdots,\, x_n)e^{-\sum_{i=1}^nf(x_i)} \d\lambda_{\phi}(x_1)\ldots\d\lambda_{\phi}(x_n)  }}\\
      &{ =\Det\left( \I-\alpha J^{\phi}_{\Lambda,\alpha}
          [e^{-f}]\right)^{-1/\alpha}.}
    \end{aligned}
  \end{multline*}
  Since $\Det\left( \I+\alpha K_{\Lambda}\right) =\Det\left( \I+\alpha
    K_{\Lambda}^\phi\right) ,$ we have:
  \begin{align*}
    \E_{\mu_{\alpha,K_{\Lambda},\lambda}} \left[ e^{- \int{f\circ \phi
          \ \d\xi}}
    \right] &=  \Det\left( \I+\alpha K^{\phi}_{\Lambda}[1-e^{-f}] \right)^{-1/\alpha}\\
    &=\E_{\mu_{\alpha,\, K_{\Lambda}^{\phi},\, \lambda_{\phi}}} \left[
      e^{- \int{f \d\xi}} \right].
  \end{align*}
  The proof is thus complete.
\end{proof}
For $\alpha=2$, Theorem \ref{thm:changementdeloi} says that the image
under $\phi$ of a Cox process is still a Cox process of parameters
$K_\Lambda^\phi$ and $\lambda_\phi$. Such a process can be constructed
as follows: Let $X$ be a centered Gaussian random field satisfying
(\ref{eq:5b}) and (\ref{eq:6}) and let $Y(x)=X(\phi^{-1}(x))$.  Then,
according to Lemma \ref{thm:operateurs}, we have: for any compact
$\Lambda$,
\begin{equation*}
  \esp{\P}{\int_\Lambda Y^2(x)\d\lambda_\phi(x)}=\trace(K_\Lambda^\phi)
\end{equation*}
and
\begin{equation*}
  \esp{\P}{Y(x)Y(y)}=K^\phi(x,\, y)=K(\phi^{-1}(x),\,
  \phi^{-1}(y)),\ \lambda_\phi\otimes \lambda_\phi, \text{ a.s..}
\end{equation*}
From Theorem \ref{thm:invariance_poisson}, by conditioning with
respect to $X,$ we also have:
\begin{align*}
  \E_{\mu_{2,\, K,\, \lambda}} \left[ e^{- \int{f\circ \phi \ \d\xi}}
  \right]&=\esp{\P}{\esp{}{\left. e^{- \int{f\circ
            \phi
            \ \d\xi}}\,\right|\, X}}\\
  &=\esp{\P}{\exp\left(-\int (1-e^{-f\circ \phi}) \, X^2\d \lambda\right)}\\
  &=\esp{\P}{\exp\left(-\int (1-e^{-f}) \, Y^2\d \lambda_\phi\right)}.
\end{align*}
Thus the two approaches (fortunately) yields the same result.

We now want to prove that $\mu_{\alpha, \, K^\phi, \, \lambda_\phi}$
is absolutely continuous with respect to $\mu_{\alpha, \, K,\,
  \lambda}$ and compute the corresponding Radon-Nikodym derivative.
For technical reasons, we need to assume that there exists a Jacobi
formula (or change of variable formula) on the measured space $(E, \,
\lambda)$. This could be done in full generality for $E$ a manifold;
for the sake of simplicity, we assume hereafter that $E$ is a domain
of some $\R^d$. We denote by $\nabla^E$ the usual gradient on
$\R^d$. We also introduce a new hypothesis.
\begin{hyp}
  \label{Assumption:rhoC1}
  We suppose that the measure $\lambda$ is absolutely continuous with
  respect to the Lebesgue measure $m$ on $E.$ We denote by $\rho$ the
  Radon-Nikodym derivative of $\lambda$ with respect to $m$.  We
  furthermore assume that $\sqrt{\rho}$ is in $H^{1,2}_{loc}(K(x,\,
  x)\d m(x))$, i.e., $\rho$ is weakly differentiable and for any
  compact $\Lambda$ in $E$, we have:
  \begin{align*}
    \infty &> 2\int_\Lambda \|\nabla^E\sqrt{\rho(x)}\|^2 K(x,\, x)\d m(x)\\
    &=  \int_\Lambda \frac{\| \nabla^E \rho(x)\|^2}{\rho(x)}K(x,\, x)\d m(x)\\
    &=\int_\Lambda \left(\frac{\| \nabla^E
        \rho(x)\|}{\rho(x)}\right)^2K(x,\, x)\d \lambda(x).
  \end{align*}
\end{hyp}
Then for any $\phi \in \Diff_0(E),$ $\lambda_\phi$ is absolutely
continuous with respect to $\lambda$ and
\begin{equation*}
  p^{\lambda}_{\phi}(x)=\frac{d\lambda_{\phi}(x)}{d\lambda(x)}
  =\frac{\rho(\phi^{-1}(x))}{\rho(x)}\Jac(\phi)(x),
\end{equation*}
where $\Jac(\phi)(x)$ is the Jacobian of $\phi$ at point $x.$
\begin{Lemma}
  \label{lem:determinant_positive}
  Assume $(E,\, K,\, \lambda)$ satisfy Hypothesis
  \ref{hyp:condition_K} and \ref{Assumption:rhoC1}. Let $(u_n, \, n\ge
  0)$ be a sequence of nonnegative real numbers such that for any
  $x\in\R,$
  \begin{equation}\label{equ:non_negative_numbers}
    \sum_{n\ge 0} \frac{u_n}{n!}\ |x|^n < +\infty.
  \end{equation}
  For any compact $\Lambda\subset E$, we have:
  \begin{equation}
    \esp{\mu_{\alpha,\, K_\Lambda, \, \lambda}}{\frac{u_{|\xi|}}{\det_{\alpha} J_{\Lambda,\alpha}(\xi)}}< + \infty.
  \end{equation}
  As a consequence, $\det_{\alpha} J_{\Lambda,\alpha}(\xi)$ is
  $\mu_{\alpha,\, K_\Lambda, \, \lambda}$ almost-surely positive.
\end{Lemma}
\begin{proof}
  According to Theorem \ref{thm:existence}, we have:
  \begin{equation*}
    j^n_{\Lambda,\alpha,K_{\Lambda}}(x_1,\cdots,\,\, x_n)=\Det (\I +\alpha K_{\Lambda})^{-1/\alpha} \det_{\alpha}  J_{\Lambda,\alpha}(x_1,\cdots,\,\, x_n),
  \end{equation*}
  hence
  \begin{multline*}
    \mathbb{E}\left[\frac{u_{|\xi|}}{\det_{\alpha}
        J_{\Lambda,\alpha}(\xi)}\right]\\
    \begin{aligned}
      &=\sum_{n=0}^{+\infty}\frac{1}{n!}\int_{\Lambda^n}\frac{u_n}{\det_{\alpha}
        J_{\Lambda,\alpha}(x_1,\cdots,\, x_n)}\ j^n_{\Lambda,\alpha,K_{\Lambda}}(x_1,\cdots,\,\, x_n)\otimes_{j=1}^n\d\lambda(x_j)\\
      &=\Det (\I +\alpha K_{\Lambda})^{-1/\alpha}
      \sum_{n=0}^{+\infty}\frac{u_n}{n!} \lambda(\Lambda)^n<+\infty,
    \end{aligned}
  \end{multline*}
  because $\lambda$ is assumed to be a Radon measure and $\Lambda$ is
  compact.
\end{proof}
\begin{Theorem}
  \label{thm:image_homeo}
  Assume $(E,\, K,\, \lambda)$ satisfy Hypothesis
  \ref{hyp:condition_K} and \ref{Assumption:rhoC1}.  Then, the measure
  $\mu_{\alpha,\, K,\, \lambda}$ is quasi-invariant with respect to
  the group $\Diff _0(E)$ and for any $\phi \in \Diff _0(E)$, we have
  then:
  \begin{equation*}
    \frac{ d\phi^*\mu_{\alpha,\,K,\,\lambda} }{
      d\mu_{\alpha,\,K,\,\lambda} }
    (\xi)=\left(\prod_{x\in\xi}p^{\lambda}_{\phi}(x)\right)
    \frac{\det_{\alpha} \, J_{\alpha}^{\phi}(\xi)}{\det_{\alpha}\, 
      J_{\alpha}(\xi)}\cdotp
  \end{equation*} 
  That is to say that for any measurable nonnegative, compactly
  supported $f$ on $E$:
  \begin{equation}\label{formule1}
    \E_{\mu_{\alpha,\, K,\, \lambda}}\left[e^{-\int f \circ \phi \d\xi}\right]=
    \E_{\mu_{\alpha,\, K,\, \lambda}}\left[e^{-\int f \d\xi }e^{\int
        ln(p^{\lambda}_{\phi})\d\xi}\ \frac{\det_{\alpha}\,
        J_{\alpha}^{\phi}(\xi)}{\det_{\alpha}\, 
        J_{\alpha}(\xi)}\right].
  \end{equation}
\end{Theorem}
\begin{proof}
  Since $f$ is compactly supported and $\phi$ belongs to $\Diff_0(E),$
  there exists a compact $\Lambda$ which contains both the support of
  $f$ and $f\circ \phi$.  According to Theorem
  \ref{thm:changementdeloi} and Lemma \ref{thm:kernel_form}, we have:
  \begin{align*}
    \E_{\mu_{\alpha,K_{\Lambda},\lambda}} \left[ e^{- \int{f\circ
          \phi\ \d\xi}} \right]
    &=\E_{\mu_{\alpha,K_{\Lambda}^{\phi},\lambda_{\phi}}} \left[ e^{- \int{f \d\xi}}    \right]\\
    &=\Det\left( \I+\alpha K_{\Lambda}^{\phi}\right)^{-1/\alpha}\left(\sum_{n=0}^{+\infty}\frac{1}{n!}A_n\right)\\
    &=\Det\left( \I+\alpha
      K_{\Lambda}\right)^{-1/\alpha}\left(\sum_{n=0}^{+\infty}\frac{1}{n!}A_n\right)
  \end{align*}
  where for any $n \in \N$, the $A_n$ are the integrals:
  \begin{multline*}
    A_n=\int_{\Lambda^n}
    {\det_{\alpha}J_{\Lambda,\alpha}^{\phi}(x_1,\, \cdots,\, x_n)
      e^{-\sum_{i=1}^n f(x_i)}\d\lambda_{\phi}(x_1)\ldots\d\lambda_{\phi}(x_n)} \\
    \begin{aligned}
      &=\int_{\Lambda^n}
      {\det_{\alpha}J^{\phi}_{\Lambda,\alpha}(x_1,\, \cdots,\, x_n)
        e^{-\sum_{i=1}^n f(x_i)} \prod_{i=1}^{n}p^{\lambda}_{\phi}(x_i)\d\lambda(x_1)\ldots \d\lambda(x_n)}\\
      &=\int_{\Lambda^n}\det_{\alpha}J_{\Lambda,\alpha}(x_1,\,
      \cdots,\, x_n) \alpha_n(x_1,\, \cdots,\,
      x_n)\d\lambda(x_1)\ldots \d\lambda(x_n) ,
    \end{aligned}
  \end{multline*}
  where
  \begin{equation*}
    \alpha_n(x_1,\, \cdots,\, x_n)=\frac{\det_{\alpha}J_{\Lambda,\alpha}^{\phi}(x_1,\, \cdots,\, x_n)}{\det_{\alpha}J_{\Lambda,\alpha}(x_1,\, \cdots,\, x_n)} e^{-\sum_{i}f(x_i)} \prod_{i=1}^{n}p^{\lambda}_{\phi}(x_i).
  \end{equation*}
  Hence according to \eqref{eq:4}, we can write:
  \begin{multline*}
    \Det\left( \I+\alpha K_{\Lambda}\right)^{-1/\alpha}\sum_{n=0}^{+\infty}\frac{1}{n!}A_n\\
    =\sum_{n=0}^{+\infty}\frac{1}{n!}
    \int_{\Lambda^n}j^n_{\Lambda,\alpha,K_{\Lambda}}(x_1,\, \cdots,\,
    x_n)\alpha_n(x_1,\, \cdots,\, x_n)\d\lambda(x_1)\ldots
    \d\lambda(x_n).
  \end{multline*}
  Thus, we have \eqref{formule1}.
\end{proof}
Should we consider Poisson process either as a $0$-DPPP or as an
$\alpha$-DPPP with the singular kernel mentioned above, we see that
the last fraction in \eqref{formule1} reduces to $1$ and we find the
well known formula of quasi-invariance for Poisson processes (see
\cite{MR99d:58179}).  In the following, we define:
\begin{equation*}
  L^{\phi}_{\mu_{\alpha,\,K,\,\lambda}}(\xi)=\left(\prod_{x \in\xi}p^{\lambda}_{\phi}(x)\right)
  \frac{\det_\alpha J_{\alpha}^{\phi}(\xi)}{\det_\alpha J_{\alpha}(\xi)}\cdotp
\end{equation*}
Then formula \eqref{formule1} can be rewritten as:
\begin{equation*}
  \E_{\mu_{\alpha,\,K,\,\lambda}} \left[ e^{- \int{f\circ \phi\
        \d\xi}}    \right]
  =\E_{\mu_{\alpha,\,K,\,\lambda}}\left[e^{-\int f \d\xi }\ L^{\phi}_{\mu_{\alpha,\,K,\,\lambda}}(\xi)\right].
\end{equation*}

\section{Integration by parts formula}\label{Integration by parts}
In this section, we prove the integration by parts formula. The proof
relies on a differentiation within (\ref{formule1}). We thus need to
put a manifold structure on $\chi.$ The tangent space $T_\xi \chi$ at
some $\xi \in \chi$ is given as $L^2(d\xi),$ i.e., the set of all maps
$V$ from $E$ to $\R$ such that:
\begin{equation*}
  \int |V(x)|^2 \d\xi(x)<\infty.
\end{equation*}
Note that if $\xi\in \chi_0$ then $T_\xi\chi$ can be identified as
$\R^{|\xi|}$ with the euclidean scalar product.

We consider $V_0(E)$ the set of all $C^{\infty}$-vector fields on $E$
with compact support. For any $v \in V_0(E)$, we construct: $\phi_t^v
: E \rightarrow E$, $t \in \R$, where the curve, for any $x \in E$
\begin{equation*}
  t \in \R \rightarrow \phi_t^v(x) 
\end{equation*}
is defined as the solution to:
\begin{equation*}
  \frac{d}{dt}\phi_t^v(x)=v(\phi^v_t(x)) \text{ and } \phi_0^v(x)=x.
\end{equation*}
Because $v \in V_0(E)$, there is no explosion and $\phi_t^v$ is
well-defined for each $t \in \R$. The mappings $\{\phi_t^v, t \in
\R\}$ form a one-parameter subgroup of diffeomorphisms with compact
support, that is to say:
\begin{itemize}
\item $\forall t \in \R, \phi_t^v \in \Diff _0(E)$.
\item $\forall t,s \in \R, \phi_t^v \circ \phi_s^v = \phi_{t+s}^v$. In
  particular, $(\phi_t^v)^{-1}=\phi_{-t}^v$.
\item For any $T>0$, there exists a compact $K$ such that
  $\phi_t^v(x)=x$ for any $x\in K^c$, for any $|t|\le T.$
\end{itemize}
In the following, we fix $v \in V_0(E)$. For any $\xi \in \chi$, we
still denote by $\phi_t^v$ the map:
\begin{align*}
  \phi_t^v\, :\, \chi & \longrightarrow \chi\\
  \xi=\sum_{x \in \xi}\delta_{x_i}&\longmapsto \sum_{x \in
    \xi}{\delta_{\phi_t^v(x)}} \in \chi.
\end{align*}
\begin{Definition}\label{def:gradient}
  A function $F :\chi \rightarrow \R$ is said to be differentiable at
  $\xi \in \chi$ whenever for any vector field $v\in V_0(E)$, the
  directional derivative along the vector field $v$
  \begin{equation*}
    \nabla_v F(\xi)=\left.\frac{d}{dt} F(\phi_t^v(\xi))\right|_{t=0}
  \end{equation*}
  is well defined.
\end{Definition}
Since $\phi_t^v$ does not change the number of atoms of $\xi$, if
$\xi$ belongs to $\chi_0$, this notion of differentiability coincides
with the usual one in $\R^{|\xi|}$ and
\begin{equation*}
  \nabla_v F(x_1,\cdots,\, x_n)=\sum_{i=1}^n \partial_i F(x_1,\cdots,\, x_n)v(x_i),
\end{equation*}
if $\xi=\{x_1,\cdots,\, x_n\}.$

In the general case, a set of test functions is defined as is :
Following the notations from \cite{MR99d:58179}, for a function
$F:\chi\rightarrow \mathbb{R}$ we say that $F \in
\mathcal{F}C^{\infty}_b(\mathcal{D},\chi)$ if:
\begin{equation*}
   F(\xi)=f\left(\int h_1 \d\xi,\, \cdots, \int h_N\d\xi  \right),
\end{equation*}
for some $N \in \N$, $h_1,\, \cdots,h_N \in
\mathcal{D}=C^{\infty}(E)$, $f \in C^{\infty}_b(\R^N)$.  Then for any
$F \in \mathcal{F}C^{\infty}_b(\mathcal{D},\chi)$, given $v \in
V_0(E)$, we have:
\begin{equation*}
  F(\phi_t^v(\xi))=f\left(\int h_1 \circ \phi_t^v\d\xi,\, \cdots, \int h_N \circ \phi_t^v\d\xi\right).
\end{equation*}
It is then clear that the directional derivative of such $F$ exists
and that:
\begin{equation*}
  \nabla_v F(\xi)=\sum_{i=1}^N \partial_i f
  \left(\int h_1 \d\xi,\, \cdots, \int h_N\d\xi \right) \int \nabla^E_v h_i\d \xi.
\end{equation*}
The gradient $\nabla F$ of a differentiable function $F$ is defined as
a map from $\chi$ into $T\chi$ such that, for any $v\in V_0(E),$
\begin{equation*}
  \int \nabla_x F(\xi) v(x)\d\xi(x)=\nabla_vF(\xi).
\end{equation*}
If $\xi\in \chi_0$ and $F$ is differentiable at $\chi$, then
\begin{equation*}
  \nabla_x F(\xi)=\sum_{i=1}^{|\xi|} \partial_i F(\{x_1,\cdots,\, x_{|\xi|}\}) \car_{\{x=x_i\}}.
\end{equation*}
If $\xi$ belongs to $\chi$, for any $F \in
\mathcal{F}C^{\infty}_b(\mathcal{D},\chi)$,
\begin{equation*}
  \nabla_x F(\xi)=\sum_{i=1}^n \partial_i f  \left(\int h_1
    \d\xi,\, \cdots, \int h_N\d\xi \right)\, \nabla^E h_i(x).
\end{equation*}
\subsection{Determinantal point processes}
\label{sec:determ-point-proc}
In what follows, $c$ and $\kappa$ are positive constant which may vary
from line to line.

In this part, we assume $\alpha=-1$ and that Hypothesis
\ref{hyp:condition_K} and \ref{Assumption:rhoC1} hold. We denote by
$\beta^{\lambda}(x)$ the logarithmic derivative of $\lambda$, given
by: for any $x$ in $E$,
\begin{equation*}
  \beta^{\lambda}(x)=\frac{\nabla\rho(x)}{\rho(x)} \text{   on }\{\rho(x)>0\},
\end{equation*}
and $\beta^{\lambda}(x)=0$ on $\{\rho(x)=0\}$.  Then, for any vector
field $v$ on $E$ with compact support, we denote by $B_v^\lambda$ the
following function on $\chi$:
\begin{align*}
  B_v^\lambda\, :\, \chi &\longrightarrow \R\\
  \xi &\longmapsto B^{\lambda}_v(\xi)=\int_E
  \left(\beta^{\lambda}(x).v(x) + \divergence (v(x))\right)\d\xi(x),
\end{align*}
where $x.y$ is the euclidean scalar product of $x$ and $y$ in $E$.  If
$\lambda=m$,
\begin{equation*}
  B_v^m(\xi)=\int_E \divergence (v(x))\d\xi(x)
\end{equation*}
and according to Theorem \ref{thm:esperance_et_variance},
\begin{align*}
  \esp{}{| B_v^m(\xi)|}&\le  \int_E |\divergence (v(x))| K(x,\, x)\d\lambda(x)\\
  &\le \|v\|_\infty \trace(K_\Lambda)<\infty,
\end{align*}
where $\Lambda$ is a compact containing the support of $v$.  As in
\cite{MR2108363}, we now define the potential energy of a finite
configuration by
\begin{align*}
  U\, :\, \chi_0 & \longrightarrow \R\\
  \xi & \longmapsto -\log \det\, J(\xi).
\end{align*}
\begin{hyp}
  \label{hyp:H_differentiable}
  The functional $U$ is differentiable at every configuration $\xi\in
  \chi_0.$ Moreover, for any $v\in V_0(E)$, there exists $c >0$ such
  that for any $\xi\in \chi_0$, we have
  \begin{equation}\label{eq:9}
    \left| \langle \nabla U(\xi),\, v\rangle_{L^2(d\xi)}\right|\le  \frac{u_{|\xi|}}{\det\, J(\xi)},
  \end{equation}
  where $(u_n=c n^{n/2},\,  n\ge 1)$ satisfy (\ref{equ:non_negative_numbers}).
\end{hyp}
\begin{Theorem}
  \label{thm:differentiability_J}
  Assume that the kernel $J$ is once differentiable with continuous
  derivative. Then, Hypothesis \ref{hyp:H_differentiable} is
  satisfied.
\end{Theorem}
\begin{proof}
  Let $\xi=\{x_1,\cdots,\, x_n\}\in \chi_0$ and let $\Lambda$ be a
  compact subset of $E$ whose interior contains $\xi$. Since $J(.,.)$
  is differentiable
  \begin{equation*}
    (y_1,\cdots,\, y_n)\longmapsto -\log\det\ (J(y_i,y_k), \, 1\le
    i,k\le n)
  \end{equation*}
  is differentiable.  The chain rule formula implies that
  \begin{equation*}
    t\longmapsto \log\det\ (J(\phi_t^v(x_i),\, \phi_t^v(x_k)), \, 1\le
    i,k\le n)
  \end{equation*}
  is differentiable and its differential is equal to
  \begin{equation*}
    \frac{1}{\det \, J(\phi_t^v(\xi))}  \trace\left( \adj  (J(\phi_t^v(x_i),\, \phi_t^v(x_k))) \left(E^v_{t}(\frac{\partial J(\xi)}{\partial x})_{t}+(\frac{\partial J(\xi)}{\partial y})_{t}E^v_{t}\right)\right),
  \end{equation*}
  where $(\frac{\partial J(\xi)}{\partial x})_{t}$ is the matrix with
  terms $\left(\frac{\partial J_{\Lambda}}{\partial
      x}(\phi^v_{t}(x_i),\phi^v_{t}(x_j))\right)_{x_i,\, x_j \in
    \xi}$, $(\frac{\partial J(\xi)}{\partial y})_{t}$ is the matrix
  with terms $\left(\frac{\partial J_{\Lambda}}{\partial
      y}(\phi^v_{t}(x_i),\phi^v_{t}(x_j))\right)_{x_i,\, x_j \in \xi
  }$, and $E^v_{t}$ is the diagonal matrix with terms
  $\left(v(\phi_{t}^v(x_i))\right)_{x_i \in \xi}$. For $t=0$, this
  reduces to
  \begin{multline*}
    \left| \langle \nabla U(\xi),\, v\rangle_{L^2(d\xi)}\right|=\\
    \frac{1}{\det \, J(\xi)} \trace\left( \adj (J(\xi))
      \left(E^v_{0}(\frac{\partial J(\xi)}{\partial
          x})_{0}+(\frac{\partial J(\xi)}{\partial
          y})_{0}E^v_{0}\right)\right).
  \end{multline*}
  Since $J$ is continuous and $\Lambda$ is compact,
  \begin{equation*}
    \| \frac{\partial J}{\partial y}(\xi)\|_{HS}\le \, |\xi|
    \|J\|_\infty \text{ and } \|E^v_0(\xi)\|_{HS}\le |\xi|^{1/2}\|v\|_\infty.
  \end{equation*}
  Hence, there exists $c$ independent of $\xi$ such that
  \begin{equation*}
    \left| \langle \nabla U(\xi),\, v\rangle_{L^2(d\xi)}\right|\le \, c\, |\xi|^{2}\frac{1}{\det \, J(\xi)} \,
    |\trace( \adj  (J(\xi)) )|.
  \end{equation*}
  From \cite[page 1021]{0635.47002}, we know that for any $n\times n$
  matrix $A$, for any $x$ and $y$ in $\R^{n}$, we have
  \begin{equation*}
    |(\adj A) x.y|\le \|y\|\|A\|_{HS}^{n-1}(n-1)^{-(n-1)/2}.
  \end{equation*}
  It follows that
  \begin{equation*}
    |\trace(\adj A)|=|\sum_{j=1}^n (\adj A) e_j.e_j|\le n \|A\|_{HS}^{n-1}(n-1)^{-(n-1)/2},
  \end{equation*}
  where $(e_j,\, j=1,\cdots,\, n)$ is the canonical basis of
  $\R^n$. Since $J$ is bounded, $\|J(\xi)\|_{HS}\le
  |\xi|\|J\|_\infty$, hence there exists $c$ independent of $\xi$ such
  that
  \begin{equation*}
    \left| \langle \nabla U(\xi),\, v\rangle_{L^2(d\xi)}\right|\le \, \frac{c}{\det \, J(\xi)} 
    |\xi|^{|\xi|/2}.
  \end{equation*}
  The proof is thus complete.
\end{proof}
\begin{Corollary}
  \label{lem:borne_sur_H}
  Assume that hypothesis \ref{hyp:H_differentiable} holds. For any
  $v\in V_0(E)$, for any $\xi\in \chi_0$, the function
  \begin{equation*}
    t\longmapsto  H_t(\xi)=\frac{\det\,  J(\phi_t^v(\xi))}{\det \, J(\xi)}
  \end{equation*}
  is differentiable and
  \begin{equation*}
    \sup_{|t|\le T} \left| \frac{d H_t(\xi)}{dt}\right|\le \frac{ u_{|\xi|}}{\det\, J(\xi)},
  \end{equation*}
  where $(u_n,\,  n\ge 0)$ satisfy (\ref{equ:non_negative_numbers}).
\end{Corollary}
\begin{proof}
  According to Hypothesis \ref{hyp:H_differentiable}, the function
  $(t\mapsto U(\phi_t^v(\xi)))$ is differentiable and
  \begin{equation}
    \label{eq:7}   \frac{d U(\phi_t^v(\xi))}{dt}=\langle \nabla U(\phi_t^v (\xi)),\
    v\rangle_{L^2(d\phi_t^v(\xi))} . 
  \end{equation}
  For any $t$, $\phi_t^v$ is a diffeomorphism hence, Theorem
  \ref{thm:image_homeo} applied to $\phi_t^v$ and $\phi_{-t}^v$
  implies that $\mu_{-1,\, K^{\phi_t^v},\lambda_{\phi_t^v}}$ and
  $\mu_{-1,\, K,\, \lambda}$ are equivalent measure. According to
  Lemma \ref{lem:determinant_positive}, for any $t$, $\det\,
  J^{\phi_t^v}(\xi)$ is $\mu_{-1,\,
    K^{\phi_t^v},\lambda_{\phi_t^v}}$-a.s. positive hence it is also
  $\mu_{-1,\, K,\, \lambda} $-a.s. positive. Since for any
  $\xi\in\chi_0$,
  \begin{equation*}
    t\mapsto \det \,
    J^{\phi_t^v}(\xi)=\exp(-U(\phi_t^v(\xi)))
  \end{equation*}
  is continuous, it follows that there exists a set of full
  $\mu_{-1,\, K,\, \lambda} $ measure on which $\det\,
  J^{\phi_t^v}(\xi)>0$ for any $|t|\le T$, for any $\xi.$ Furthermore,
  \begin{align*}
    \frac{d H_t(\xi)}{dt} = -\frac{\det\, J(\phi_t^v(\xi))}{\det\,
      J(\xi)}\ \frac{d U(\phi_t^v(\xi))}{dt}.
  \end{align*}
  In view of (\ref{eq:7}) and of Hypothesis
  \ref{hyp:H_differentiable}, this means that
  \begin{align*}
    \sup_{|t|\le T} \left| \frac{d H_t(\xi)}{dt}\right|&\le
    \frac{\det\, J(\phi_t^v(\xi))}{\det\,
      J(\xi)}\ \frac{u_{|\xi|}}{\det\, J(\phi_t^v(\xi))}\\
    &= \frac{ u_{|\xi|}}{\det\, J(\xi)},
  \end{align*}
  since $\phi_t^v(\xi)$ has the number of atoms as $\xi.$
\end{proof}
\begin{Lemma}
  \label{lem:1}
  Assume that $\lambda=m$ and set
  \begin{equation*}
    P_t(\xi)=\prod_{x\in \xi} p_{\phi_t^v}(x)=\prod_{x\in \xi} \Jac \phi_t^v(x).
  \end{equation*}
  For any $v\in \Diff_0(E)$, for any configuration $\xi \in \chi$, $P$
  is differentiable with respect to $t$ and we have
  \begin{equation*}
    \frac{d \log P_t}{dt}(\xi)  =\int \left(  \divergence v -\int_0^t
      \nabla^E \divergence v \circ \eta_{r,t} .\ v(\eta_{r,t}) \d r \right) d\xi,
  \end{equation*}
  where for any $r\le t,\ x\rightarrow\eta_{r,t}(x)$ is the
  diffeomorphism of $E$ which satisfies:
  \begin{equation*}
    \eta_{r,t}(x)=x-\int_r^t v(\eta_{s,t}(x))\d s.
  \end{equation*}
  In particular for $t=0$, we have:
  \begin{equation}
    \frac{d}{dt}\left.\left(\prod_{x \in\xi}p^{\lambda}_{\phi_t^v}(x)\right)\right|_{t=0}=B^{m}_v(\xi).
  \end{equation}
  Moreover, there exists $c>0$ and $\kappa >0$ such that for any $\xi
  \in \chi_0$,
  \begin{equation}\label{eq:14}
    \sup_{t\le T} \left | \frac{d  P_t}{dt}(\xi)\right|\le c e^{\kappa |\xi|}.
  \end{equation}
\end{Lemma}
\begin{proof}
  Introduce, for any $s\le t,\ x\longmapsto\eta_{s,t}(x)$, the
  diffeomorphism of $E$ which satisfies:
  \begin{equation*}
    \eta_{s,t}(x)=x-\int_s^t v(\eta_{r,t}(x))\d r.
  \end{equation*}
  As a comparison, we remind that $\phi_t^v(x)=x+\int_0^t
  v(\phi_s^v(x))\d s$. It is well-known that the diffeomorphism
  $x\longmapsto \eta_{0,t}(x)$ is the inverse of $x\longmapsto
  \phi^v_t(x)$. Then using \cite{ustunel2000}, we have:
  \begin{equation}\label{eq:11}
    \Jac\phi_t^v(x)=\frac{d(\phi_t^v)^*m(x)}{dm(x)}=
    \exp\left(\int_0^t\divergence \ v\circ\eta_{r,t}(x)\d r\right),
  \end{equation}
  and:
  \begin{equation*}
    \prod_{x \in \xi}\Jac\phi_t^v(x) =\exp\left(\sum_{x \in
        \xi}\int_0^t\divergence \ v\circ\eta_{r,t}(x)\d r\right).
  \end{equation*}
  Hence, we have:
  \begin{align*}
    \sum_{x\in \xi} \frac{d}{dt} \log \Jac \phi_{t}^v(x)&=\sum_{x\in
      \xi}  \frac{d}{dt} \int_0^t\divergence \ v\circ\eta_{r,t}(x)\d r\\
    & =\sum_{x \in\xi}\divergence v(x) -\int_0^t\nabla^E\divergence v
    \circ \eta_{r,t}(x). v(\eta_{r,t}(x)) \d r.
  \end{align*}
  The first and second point follow easily. Now, $v$ is assumed to
  have bounded derivatives of any order, hence for any $\xi\in
  \chi_0$,
  \begin{equation}\label{eq:12}
    \left| \frac{d\log P_t}{dt}(\xi)\right| \le c |\xi|,
  \end{equation}
  where $c$ does not depend neither from $t$ nor $\xi$. According
  to~(\ref{eq:11}), there exists $\kappa >0$ such that for any $\xi\in
  \chi_0$, we have:
  \begin{equation}\label{eq:13}
    |P_t(\xi)|\le \exp(\kappa |\xi|).
  \end{equation}
  Thus, combining~(\ref{eq:12}) and~(\ref{eq:13}), we
  get~(\ref{eq:14}).
\end{proof}
We are now in position to prove the main result of this section.
\begin{Theorem}
  \label{thm:thm_principal}
  Assume $(E,\, K,\, \lambda)$ satisfy Hypothesis
  \ref{hyp:condition_K}, \ref{Assumption:rhoC1} and
  \ref{hyp:H_differentiable}, let $\alpha=-1$. Let $F$ and $G$ belong
  to $\FC.$ For any compact $\Lambda$, we have:
  \begin{multline}
    \label{eq:10}
    \int_{\chi_\Lambda} \nabla_v F(\xi) G(\xi)\d\mu_{-1,K_\Lambda,\lambda}(\xi) =-\int_{\chi_\Lambda} F(\xi)\nabla_v G(\xi)\d\mu_{-1,K_\Lambda,\lambda}(\xi)\\
    +\int_{\chi_\Lambda} F(\xi)G(\xi)\left(B_v^\lambda
      (\xi)+\nabla_vU(\xi)\right)\d\mu_{-1,K_\Lambda,\lambda}(\xi).
  \end{multline}
\end{Theorem}
\begin{proof}
  In view of Lemma \ref{lem:chgt_mesure_operateur}, we can replace $J$
  by $J[\rho]$ and assume $\lambda=m$, i.e., $\lambda$ is the Lebesgue
  measure. Note that
  \begin{equation*}
    B^m_v(\xi)=\int \divergence v(x)\d\xi(x).
  \end{equation*}
  Let $\Lambda$ be a fixed compact set in $E,$ remember that
  $\chi_\Lambda\subset \chi_0$. Let $M$ be an integer and
  $\chi^M=\{\xi\in \chi_0,\ |\xi|\le M\}.$ It is crucial to note that
  $\chi^M$ is invariant by any $\phi\in \Diff_0(E)$.  On the one hand,
  by dominated convergence, we have:
  \begin{multline*}
    \frac{d}{dt}\left.\left(\int_{\chi^M}F(\phi_t^v(\xi))G(\xi)\d\mu_{-1,K_\Lambda[\rho],m}(\xi)\right)\right|_{t=0}\\
    \begin{aligned}
      &=\int_{\chi^M}\left.\frac{d}{dt}\left(F(\phi_t^v(\xi))\right)\right|_{t=0}G(\xi)\d\mu_{-1,K_\Lambda[\rho],m}(\xi)\\
      &=\int_{\chi^M}\nabla_v
      F(\xi)G(\xi)\d\mu_{-1,K_\Lambda[\rho],m}(\xi).
    \end{aligned}
  \end{multline*}
  On the other hand, we know from (\ref{formule1}) that
  \begin{multline}\label{eq:15}
    \int_{\chi^M}F(\phi_t^v(\xi))G(\xi)\d\mu_{-1,K_\Lambda[\rho],m}(\xi)\\
    \begin{aligned}
      &=\int_{\chi_\Lambda}
      F(\phi_t^v(\xi))G(\xi)\car_{\{|\xi|\le M\}}\d\mu_{-1,K_\Lambda[\rho],m}(\xi)\\
      &=\int_{\chi_\Lambda} F(\xi)G(\phi_{-t}^v(
      \xi))\car_{\{|\phi_{-t}^v( \xi)|\le M\}}\d
      \mu_{-1,K^{\phi_t^v}_\Lambda[\rho],m_{\phi_t^v}}(\xi)\\
      &=\int_{\chi_\Lambda} F(\xi)G(\phi_{-t}^v( \xi))\car_{\{|
        \xi|\le M\}}L_{-1,\, K[\rho],\,
        \lambda}^{\phi_t^v}(\xi)\d\mu_{-1,K_\Lambda[\rho],m}(\xi).
    \end{aligned}
  \end{multline}
  According to Corollary \ref{lem:borne_sur_H} and Lemma \ref{lem:1},
  the function $(t\mapsto L_{-1,\, K[\rho],\, \lambda}^{\phi_t^v}
  (\xi))$ is differentiable and there exists $c$ such that:
  \begin{equation*}
    \sup_{t\le T}\left |\frac{d L_{-1,\, K[\rho],\, \lambda}^{\phi_t^v} }{dt}(\xi)\right|\le
    \frac{ u_{|\xi|}}{\det\, J(\xi)},
  \end{equation*}
  where $(u_n, n\ge 0)$ satisfy (\ref{equ:non_negative_numbers}).\\
  Lemma \ref{lem:determinant_positive} implies that the
  right-hand-side of the last inequality is integrable with respect to
  $\mu_{-1,\, K_\Lambda, \, \lambda},$ thus, we can differentiate
  inside the expectations in (\ref{eq:15}) and we obtain:
  \begin{multline*}
    \int_{\chi_\Lambda} \nabla_v F(\xi) G(\xi) \car_{\{| \xi|\le
      M\}}\d\mu_{-1,K_\Lambda,m}(\xi) \\
    =\int_{\chi_\Lambda} F(\xi)\left(-\nabla_v G(\xi) +G(\xi)
      \left(B_v^m (\xi)+\nabla_vU(\xi)\right) \right)\car_{\{| \xi|\le
      M\}}\d\mu_{-1,K_\Lambda,m}(\xi).
  \end{multline*}
  According to Hypothesis \ref{hyp:H_differentiable} and Lemma
  \ref{lem:determinant_positive}, by dominated convergence, we have:
  \begin{multline*}
    \int_{\chi_\Lambda} \nabla_v F(\xi) G(\xi) \d\mu_{-1,K_\Lambda,m}(\xi) \\
    =\int_{\chi_\Lambda} F(\xi)\left(-\nabla_v G(\xi) +G(\xi)
      \left(B_v^m (\xi)+\nabla_vU(\xi)\right)
    \right)\d\mu_{-1,K_\Lambda,m}(\xi).
  \end{multline*}
  Now, we remark that
  \begin{align*}
    \nabla_v U[\rho](\xi)&=\nabla_v \log \det J[\rho](\xi)\\
    &=\nabla_v \log \left (\prod_{x\in \xi}\rho(x)\det J(\xi)\right)\\
    &=\nabla_v \int \log \rho(x)\d\xi(x)+\nabla_v U(\xi)\\
    &=\int \frac{\nabla^E\rho(x)}{\rho(x)}.v(x)\d\xi(x) +\nabla_v
    U(\xi).
  \end{align*}
  Moreover, we have
  \begin{equation*}
    B_v^m (\xi)+\int_\Lambda \frac{\nabla^E\rho(x)}{\rho(x)}.v(x)\d\xi(x)=B^\lambda_v(\xi),
  \end{equation*}
  and in view of Theorem \ref{thm:esperance_et_variance},
  \begin{multline*}
    \esp{}{\left|\int_\Lambda
        \frac{\nabla^E\rho(x)}{\rho(x)}.v(x)\d\xi(x)\right|}^2\\
    \begin{aligned}
      & \le \esp{}{\int_\Lambda
        \left(\frac{\|\nabla^E\rho(x)\|}{\rho(x)}\right)^2\d\xi(x)}\esp{}{\int_\Lambda
        |v(x)|^2\d\xi(x)}
      \\
      &\le \|v\|_\infty^2 \trace(K_\Lambda)\ \int_\Lambda
      \left(\frac{\|\nabla^E\rho(x)\|}{\rho(x)}\right)^2K(x,\,
      x)\rho(x)\d m(x).
    \end{aligned}
  \end{multline*}
  Then, Hypothesis \ref{Assumption:rhoC1} implies that $B^\lambda_v$
  is integrable and we get \eqref{eq:10} in the general case.
\end{proof}
\subsection{$\alpha$-determinantal point processes}
\label{sec:alpha-determ-point}
We now prove the integration by parts formula for
$\alpha$-determinantal point processes where $\alpha=-1/s$ for $s$
integer greater than $2$. In principle, we could follow the previous
lines of proof modifying the definition of $U$ as
\begin{equation*}
  U(\xi)=-\log\det_\alpha\, J_\alpha(\xi)
\end{equation*}
and assuming that Hypothesis \ref{hyp:H_differentiable} is still
valid. Unfortunately, there is no (simple) analog of Theorem
\ref{thm:differentiability_J} since there is no rule to differentiate
an $\alpha$-determinant and control its derivative.

We already saw that such an $\alpha$-DPPP can be obtained as the
superposition of $s$ determinantal processes of kernel $K/s$.

Let $(E_1,\, \lambda_1, \, K_1),\cdots, \, (E_s, \,\lambda_s,\, K_s)$
be $s$ Polish spaces each equipped with a Radon measure and $s$ linear
operators satisfying Hypothesis \ref{hyp:condition_K} on their
respective space. We set
\begin{equation*}
  E=\cup_{i=1}^s \{i\}\times E_i,
\end{equation*}
that is to say $E$ is the disjoint union of the $E_i$'s, often denoted
as $\sqcup_{i=1}^s E_i$. An element of $E$ is thus a couple $(i,\, x)$
where $x$ belongs to $E_i$ for any $i\in \{1,\cdots,\, s\}.$ On the
Polish space $E$, we put the measure $\lambda$ defined by
\begin{equation*}
  \int_E f(i,\, x)\d\lambda(i,\ x)=\int_E f(i,\, x)\d\lambda_i(x).
\end{equation*}
We also define $K$ as
\begin{equation*}
  Kf(i,\, x)=\int_{E_i} K_i(x,\, y)f(y)\d\lambda_i(y).
\end{equation*}
A compact set in $E$ is of the form $\Lambda=\cup_{i=1}^s \{i\}\times
\Lambda_i$ where $\Lambda_i$ is a compact set of $E_i$ hence
\begin{equation*}
  K_\Lambda f(i,\, x)=\int_{\Lambda_i} K_i(x,\, y)f(y)\d\lambda_i(y).
\end{equation*}
This means that $K$ is a kernel operator the kernel of which is given
by:
\begin{equation}\label{eq:17}
  K((i,\, x),\, (j,\, y))=K_i(x,\, y) \car_{\{i=j\}}.
\end{equation}
In particular, for $\xi=((i_l,\, x_l),\ l=1,\cdots,\, n),$ we have
\begin{equation*}
  \det \,K(\xi)=\prod_{j=1}^s \det \, K(\xi_j)
\end{equation*}
where $\xi_j=\{x,\ (j,\, x)\in \xi\}.$

It is straightforward that $K$ is symmetric and locally of trace class. Moreover,
its spectrum is equal to the union of the spectra of the $K_i$'s. For,
if $\psi$ is such that $K\psi=\alpha\psi$ then $\psi(i,.)$ is an
eigenvector of $K_i$ and thus $\alpha$ belongs to the spectrum of
$K_i$. In the reverse direction, if $\psi$ is an eigenvector of $K_i$
associated to the eigenvalue $\alpha$ then the function
\begin{equation*}
  f(j,\, x)=\psi(x)\car_{\{i=j\}}
\end{equation*}
is square integrable with respect to $\lambda$ and is an eigenvector
of $K$ for the eigenvalue $\alpha.$ If we assume furthermore that each
of the $E_i$'s is a subset of $\R^d$, we can define the gradient on
$E$ as
\begin{equation*}
  \nabla^E f(i,\, x)=\nabla^{E_i}f(i,\, x).
\end{equation*}
Now $\chi_E$ is the set of locally finite point measures of the form
\begin{equation*}
  \xi=\sum_j \delta_{(i_j,\, x_j)}.
\end{equation*}
With these notations, it is clear that Hypothesis
\ref{hyp:condition_K}, \ref{Assumption:rhoC1} and
\ref{hyp:H_differentiable} are satisfied provided they are satisfied
for each index $i$. Thus (\ref{eq:10}) is satisfied.

Now take $E_1=\ldots=E_s$, $\lambda_1=\ldots=\lambda_s$ and
$K_1=\ldots=K_s$. We introduce the map $\Theta$ defined as:
\begin{align*}
  \Theta \, :\, E&\longrightarrow E_1\\
  (i,\, x) & \longmapsto x.
\end{align*}
Consistently with earlier defined notations, we still denote by
$\Theta$ the map
\begin{align*}
  \Theta \, :\, \chi_E & \longrightarrow \chi_{E_1}\\
  \xi & \longmapsto \sum_{(j,\, x)\in \xi} \delta_{x}.
\end{align*}
Then, according to what has been said above, $\mu_{-1/s,\, sK_1,\,
  \lambda_1}$ is the image measure of $\mu_{-1, \, K,\, \lambda}$ by
the map $\Theta.$ Set
$$\xi_n=\sum_{(i,\, x)\in\xi} \delta_x\car_{\{i=n\}}.$$ 
The reciprocal problem, interesting in its own sake and useful for the
sequel, is to determine the conditional distribution of $\xi_n$ given
$\Theta\xi.$
\begin{Theorem}
  \label{thm:conditional_distribution}
  Let $s$ be an integer strictly greater than $1$, for $F$
  non-negative or bounded, for any $\Lambda$ compact subset of $E,$
  \begin{equation}\label{eq:16}
    \esp{}{F(\xi_1)\, |\, \Theta\xi}=\sum_{\eta\subset \Theta\xi} F(\eta)\times 
    \binom{|\Theta\xi|}{|\eta|} \frac{j_{\beta,\,
        (s-1)K_{1,\Lambda}}(\Theta\xi\backslash \eta)\  j_{-1,\,
        K_{1,\Lambda}}(\eta)}{j_{\alpha, sK_{1,\Lambda}}(\Theta\xi)},
  \end{equation}
  where $\beta=-1/(s-1).$ 

  Note that (\ref{eq:16}) also holds for $s=1$ with the convention
  that $j_{\beta,\,0}(\eta)=0$ for $\eta\neq \emptyset$ and
  $j_{\beta,\,0}(\emptyset)=1$, which is analog to the usual
  convention $0^0=1$.
\end{Theorem}
\begin{proof}
  Let $\zeta=\xi_2\cup\ldots\cup\xi_s,$ we known that $\zeta$ is
  distributed as $\mu_{-\beta, \, - K_1/\beta,\, \lambda_1}.$ Consider
  $\Xi,$ the map
  \begin{align*}
    \Xi \, :\, \chi_{E_1}\times \chi_{E_1} & \longrightarrow
    \chi_{E_1}\times \chi_{E_1}\\
    (\eta_1,\, \eta_2) &\longmapsto (\eta_1, \, \eta_1\cup \eta_2).
  \end{align*}
  By construction, the joint distribution of $\Xi(\xi_1,\, \zeta)$ is
  the same as the distribution of $(\xi_1,\, \Theta\xi)$. For any
  $\eta\subset \Theta\xi\in \chi_0$, we set:
  \begin{equation*}
    R(\eta,\, \Theta\xi)=\binom{|\Theta\xi|}{|\eta|} \frac{j_{\beta,\,
        (s-1)K_{1,\Lambda}}(\Theta\xi\backslash \eta)\   j_{-1,\,
        K_{1,\Lambda}}(\eta)}{j_{\alpha, sK_{1,\Lambda}}(\Theta\xi)}\cdotp
  \end{equation*}
  Hence, we have
  \begin{multline*}
    \esp{}{F(\xi_1)G(\Theta\xi)} =\esp{}{(F\otimes G)\circ
      \Xi(\xi_1,\,
      \zeta)}\\
    = \sum_{j=0}^\infty\sum_{k=0}^\infty
    \frac{1}{j!}\frac{1}{k!}\int_{\Lambda^j\times \Lambda^k}
    \kern-5pt  F(\{x_1,\cdots, x_j\})G(\{x_1,\cdots, x_j\}\cup\{y_1,\cdots,y_k\})\\
    \shoveright{ \times j_{-1,\, K_{1,\Lambda}}(x_1,\cdots, \, x_j)\
      j_{\beta,\,  (s-1)K_{1,\Lambda}}(y_1,\cdots,\,y_k)\d\lambda(x_1)\ldots\d\lambda(y_k)}\\
    = \sum_{j=0}^\infty\sum_{k=0}^\infty
    \frac{1}{(k+j)!}\int_{\Lambda^j\times \Lambda^k}
    \kern-5pt F(\{x_1,\cdots, x_j\})(GR)(\{x_1,\cdots, x_j\}\cup\{y_1,\cdots,y_k\})\\
    \shoveright{ \times j_{\alpha, sK_{1,\Lambda}}(x_1,\cdots, x_j
      ,\,y_1,\cdots,\,y_k)\d\lambda(x_1)\ldots\d\lambda(y_k)}\\
    \shoveleft{= \sum_{m=0}^\infty \frac{1}{m!}\int_{\Lambda^m}
      \left(\sum_{j\le m}
        F(\{x_1,\cdots,\, x_j\})R(\{x_1,\cdots,\, x_j\},\, \{x_1,\cdots,\, x_m\}) \right)}\\
    \shoveright{\times G(\{x_1,\cdots,\, x_m\})\, j_{\alpha,
        sK_{1,\Lambda}}(x_1,\cdots,\, x_m)
      \d\lambda(x_1)\ldots\d\lambda(x_m)}\\
    =\int_{\chi_{E_1}} \left(\sum_{\eta\subset \omega} F(\eta)
      R(\eta,\ \omega) \right)\, G(\omega)\d\mu_{\alpha,\, sK_1,\,
      \lambda_1}(\omega).
  \end{multline*}
  The proof is thus complete.
\end{proof}
This formula can be understood by looking at the extreme case of
Poisson process. Assume that $\Theta\xi$ is distributed according to a
Poisson process of intensity $\lambda\d m$. Then, $\xi_1$ is a Poisson
process of intensity $s^{-1}\lambda \d m$ and $\zeta$ also is a
Poisson process of intensity $(1-s^{-1})\lambda \d m$. The couple
$(\xi_1,\, \Theta\xi)$ can then be constructed by random thinning of
$\Theta\xi$: Keep each point of $\Theta\xi$ independently of the
others, with probability $1/s$; the remaining points will be
distributed as $\xi_1$. The conditional expectation of a functional
$F(\xi_1)$ given $\Theta\xi$ is then the sum of the values of $F$
taken for each realization of a thinning multiplied by the probability
of each thinned configuration. Since $|\Theta\xi|$ is assumed to be
known, the atoms of $\Theta\xi$ are independent and identically
dispatched along $E$, hence the probability to obtain a specific
configuration is binomially distributed of parameters $|\Theta\xi|$
and $1/s.$ This means that
\begin{equation*}
  \esp{}{F(\xi_1)\, |\, \Theta\xi}=\sum_{\eta\subset \Theta\xi} F(\eta)\times 
  \binom{|\Theta\xi|}{|\eta|} \left(\frac{1}{s}\right)^{|\eta|}\left(1-\frac 1s\right)^{|\Theta\xi|-|\eta|}.
\end{equation*}
This corresponds to (\ref{eq:16}) for $\alpha=0$.  As a consequence,
(\ref{eq:16}) can be read as a generalization of this procedure where
the points cannot be drawn independently and with equal probability
because of the correlation structure.

For $h$ any map from $E_1$ into $E_1$, we define $h^\sqcup$ by
\begin{align*}
  h^\sqcup \,: \, E&\longrightarrow E\\
  (i,\, x) &\longmapsto (i,\, h(x)).
\end{align*}
With this notation at hand, for $v$ in $V_0(E_1)$, $(\phi_t^v)^\sqcup$
is the solution of the equations:
\begin{equation*}
  \d (\phi_t^v)^\sqcup(i,\, x)=v^\sqcup((\phi_t^v)^\sqcup(i,\, x)), \
  1\le i \le m.
\end{equation*}
Note that we only consider a restricted set of perturbations of
configurations in the sense that we move atoms on each ``layers''
without ``crossing'': By the action of $(\phi_t^v)^\sqcup$, an atom of
the form $(i,\, x)$ is moved into an atom of the form $(i,\, y)$,
leaving its first coordinate untouched.
\begin{Theorem}
  \label{thm:ipp_alpha_determinantal}
  Assume that $(E_1,\, K_1, \lambda_1)$ satisfy Hypothesis
  \ref{hyp:condition_K}, \ref{Assumption:rhoC1} and
  \ref{hyp:H_differentiable}. Let $s=-1/\alpha$ be an integer greater
  than $1$. For $F$ and $G$ cylindrical functions, for $v\in
  V_0(E_1)$, we have:
  \begin{multline*}
    \int_{\chi_\Lambda} \nabla_{v} F(\omega) G(\omega)\d\mu_{\alpha,\,
      sK_{1,\Lambda},\lambda_1}(\omega) =-\int_{\chi_\Lambda}
    F(\omega)\nabla_{v} G(\omega)\d\mu_{\alpha,\, K_{1,\,\Lambda},\lambda_1}(\omega)\\
    +\frac 1{|\alpha|}\int_{\chi_\Lambda}
    F(\omega)G(\omega)\left(\sum_{\eta\subset\omega}
      (B^{\lambda_1}_v(\eta)+\nabla_vU(\eta))R(\eta,\, \omega)
    \right)\d\mu_{\alpha,\,sK_{1,\Lambda},\lambda_1}(\omega).
  \end{multline*}
\end{Theorem}
\begin{proof}
  We first apply (\ref{eq:10}) to the process $\xi=(\xi_1,\cdots,\,
  \xi_s)$. Remember that $\Theta \xi$ is equal to $
  \xi_1\cup\ldots\cup\xi_s$. A cylindrical function of $\Theta \xi$ is
  a function of the form:
  \begin{equation*}
    H(\Theta \xi)=f(\int h_1\d \Theta \xi,\cdots,\, \int h_N\d \Theta \xi)
  \end{equation*}
  where $h_1,\, \cdots,h_N \in \mathcal{D}=C^{\infty}(E_1)$, $f \in
  C^{\infty}_b(\R^N)$. Such a functional can be written as
  \begin{math}
    F\circ \Theta(\xi)
  \end{math}
  where $F$ is a cylindrical function of $\xi.$ Moreover, for $v\in
  V_0(E_1)$,
  \begin{align}
    \nabla_v H(\Theta \xi)&=\lim_{t\to
      0}\frac{1}{t}\left(H(\phi_t^v(\Theta \xi )-H(\Theta \xi)\right)\notag\\
    &=\lim_{t\to
      0}\frac{1}{t}\left(F(\Theta(\phi_t^v)^\sqcup(\xi)-F(\Theta\xi)\right)\notag\\
    &=\nabla_{v^\sqcup}F(\Theta \xi).\label{eq:18}
  \end{align}
  In view of (\ref{eq:17}),
  \begin{equation}\label{eq:20}
    U(\xi)=-\log\det\, J(\xi_1,\cdots,\,\xi_s)=\sum_{j=1}^s U(\xi_j).
  \end{equation}
  Analyzing the proof of (\ref{eq:10}), we see that the intrinsic
  definition of $B^\lambda_v$ is
  \begin{equation*}
    B^\lambda_v(\xi)=\int \divergence_\lambda(v)\d\xi
  \end{equation*}
  where
  \begin{equation*}
    \divergence_\lambda(v)(x)=\left.\frac{d}{dt}\left(\frac{d(\phi_t^v) ^*\lambda}{\d\lambda}(x)\right)\right|_{t=0}.  
  \end{equation*}
  In view of (\ref{eq:18}), we only need to consider flows on $E$
  associated to vector fields of the form $v^\sqcup$ for $v\in
  V_0(E_1).$ Hence,
  \begin{equation}\label{eq:19}
    B^\lambda_{v^\sqcup}(\xi)=\sum_{j=1}^s B_v^{\lambda_j}(\xi_j).
  \end{equation}
  It follows from the previous considerations that:
  \begin{multline*}
    \int_{\chi_{\Lambda^\sqcup}} \nabla_{v^\sqcup} F(\Theta\xi) G(\Theta\xi)\d\mu_{-1,K_\Lambda,\lambda}(\xi) =-\int_{\chi_{\Lambda^\sqcup}} F(\Theta\xi)\nabla_{v^\sqcup} G(\Theta\xi)\d\mu_{-1,K_\Lambda,\lambda}(\xi)\\
    +\int_{\chi_{\Lambda^\sqcup}}
    F(\Theta\xi)G(\Theta\xi)\left(B_{v^\sqcup}^\lambda
      (\xi)+\nabla_{v^\sqcup}U(\xi)\right)\d\mu_{-1,K_\Lambda,\lambda}(\xi)
  \end{multline*}
  where ${\Lambda^\sqcup}=\cup_{j=1}^s \{i\}\times \Lambda.$ Since the
  $\xi_j$'s are independent and identically distributed, according to
  (\ref{eq:20}) and (\ref{eq:19}), we have
  \begin{align*}
    \esp{}{\left.B_{v^\sqcup}^\lambda (\xi)+\nabla_{v^\sqcup}U(\xi)\,
      \right|\, \Theta\xi}&=s
    \esp{}{\left.B^{\lambda_1}_v(\xi_1)+\nabla_vU(\xi_1)\,
      \right |\, \Theta\xi}\\
    &=-\frac 1\alpha\sum_{\eta\subset\Theta\xi}
    (B^{\lambda_1}_v(\eta)+\nabla_vU(\eta))R(\eta,\, \Theta\xi).
  \end{align*}
  Thus, we obtain:
  \begin{multline*}
    \int_{\chi_\Lambda} \nabla_{v} F(\omega) G(\omega)\d\mu_{\alpha,\,
      sK_{1,\Lambda},\lambda_1}(\omega) =-\int_{\chi_\Lambda}
    F(\omega)\nabla_{v} G(\omega)\d\mu_{\alpha,\, K_{1,\,\Lambda},\lambda_1}(\omega)\\
    -\frac 1\alpha\int_{\chi_\Lambda}
    F(\omega)G(\omega)\left(\sum_{\eta\subset\omega}
      (B^{\lambda_1}_v(\eta)+\nabla_vU(\eta))R(\eta,\, \omega)
    \right)\d\mu_{\alpha,\,sK_{1,\Lambda},\lambda_1}(\omega).
  \end{multline*}
  The proof is thus complete.
\end{proof}
\subsection{$\alpha$-permanental point processes}
\label{sec:alpha-perm-point}
For permanental point processes, we begin with the situation where
$\alpha=1$. In this case,
\begin{equation*}
  j_{1,\, K_\Lambda,\,
    \lambda}(\{x_1,\cdots,\, x_n\})=\Det(I+K_\Lambda)^{-1} \text{per}
  (J(x_i,\, x_j), \, 1\le i,j\le n).
\end{equation*}
We aim to follow the lines of proof of Theorem
\ref{thm:thm_principal}, for, we need some preliminary considerations.

For any integer $n$, let $D[n]$ be the set of partitions of
$\{1,\cdots,\, n\}.$ The cardinal of $D[n]$ is known to be the $n$-th
Bell number (see \cite{MR1434477}), denoted by $\BB_n$ and which can
be computed by their exponential generating function: for any real
$x$,
\begin{equation}
  \label{eq:23}
  \sum_{n=0}^\infty \BB_n \frac{x^n}{n!}=e^{e^x}-1.
\end{equation}
For an $n\times n$ matrix $A=(a_{ij}, \, 1\le i,\, j\le n)$ and for
$\tau$ a subset of $\{1,\cdots,\,n \}$, we denote by $A[\tau]$ the
matrix $(a_{ij},\, i\in \tau,\, j\in \tau)$. For a partition $\sigma$
of $\{1,\cdots,\, n\}$, $\iota(\sigma)$ is the number of non-empty
parts of $\sigma$. This means that $\sigma=(\tau_1,\cdots,\,
\tau_{\iota(\sigma)}),$ where the $\tau_i$'s are disjoint subsets of
$\{1,\cdots,\, n\}$ whose union is exactly $\{1,\cdots,\, n\}$. Then,
we set
\begin{equation*}
  \det\, A[\sigma]= \prod_{j=1}^{\iota(\sigma)} \det \, J[\tau_j].
\end{equation*}
It is proved in \cite[Corollary 1.7]{MR1433237} that
\begin{equation}\label{eq:24}
  \per A=\sum_{\sigma\in D[n]} (-1)^{n+\iota(\sigma)}\, \det\ A[\sigma].
\end{equation}
We slightly change the definition of the potential energy of a finite
configuration as
\begin{align*}
  U\, :\, \chi_0 & \longrightarrow \R\\
  \xi & \longmapsto -\log \per\, J(\xi).
\end{align*}
A new hypothesis then arises:
\begin{hyp}
  \label{hyp:U_differentiable}
  The functional $U$ is differentiable at every configuration $\xi\in
  \chi_0.$ Moreover, for any $v\in V_0(E)$, there exists $(u_n,\, n\ge
  1)$ a sequence of nonnegative real as in Lemma
  \ref{lem:determinant_positive} such that for any $\xi\in \chi_0$, we
  have
  \begin{equation}\label{eq:9}
    \left| \langle \nabla U(\xi),\, v\rangle_{L^2(d\xi)}\right|\le
    \frac{u_{|\xi|}}{\per\, J(\xi)}\cdotp
  \end{equation}
\end{hyp}
An analog of Theorem \ref{thm:differentiability_J} now becomes
\begin{Theorem}
  \label{thm:differentiability_J_per}
  Assume that $K$ is of finite rank $N$ and that the kernel $J$ is
  once differentiable with continuous derivative. Then, Hypothesis
  \ref{hyp:U_differentiable} is satisfied.
\end{Theorem}
\begin{proof}
  Since $K$ is of finite rank $N$ there are at most $N$ points in any
  configuration.  It is clear from (\ref{eq:24}) that $(t\mapsto
  U(\phi_t^v(\xi)))$ is differentiable.  Since $|\det J(\xi)[\tau]|\le
  c|\tau|^{|\tau/2|}$ where $|\tau|$ is the cardinal of $\tau\in
  D[|\xi|]$, we get
  \begin{equation*}
    \left|\langle \nabla U(\xi),\, v\rangle_{L^2(d\xi)}\right|\le c \,
    \frac{\BB_{|\xi|} |\xi|^{|\xi|/2}}{\per \, J(\xi)}\car_{\{|\xi|\le N\}}.
  \end{equation*}
  Hence the result.
\end{proof}
\begin{Remark}
  The finite rank condition is rather restrictive but the sequence
  $(\BB_n n^{n/2},\, n\ge 1)$ has not a finite exponential generating
  function thus we can't avoid it. In order to circumvent this
  difficulty one would have to improve known upper-bounds on
  permanents.
\end{Remark}
We can then state the main result for this subsection.
\begin{Theorem}
  \label{thm:ipp_permanental}
  Assume that $(E,\, K, \, \lambda)$ satisfy Hypothesis
  \ref{hyp:condition_K}, \ref{Assumption:rhoC1} and
  \ref{hyp:U_differentiable}.  Let $F$ and $G$ belong to $\FC.$ For
  any compact $\Lambda$, we have:
  \begin{multline*}
    \int_{\chi_\Lambda} \nabla_v F(\xi) G(\xi)\d\mu_{1,K_\Lambda,\lambda}(\xi) =-\int_{\chi_\Lambda} F(\xi)\nabla_v G(\xi)\d\mu_{1,K_\Lambda,\lambda}(\xi)\\
    +\int_{\chi_\Lambda} F(\xi)G(\xi)\left(B_v^\lambda
      (\xi)+\nabla_vU(\xi)\right)\d\mu_{1,K_\Lambda,\lambda}(\xi).
  \end{multline*}
\end{Theorem}
\begin{proof}
  Same as the proof of Theorem \ref{thm:thm_principal}.
\end{proof}
Now then, we can work as in Subsection \ref{sec:alpha-determ-point}
and we obtain the integration by parts formula for
$\alpha$-permanental point processes.
\begin{Corollary}
  Assume that $(E_1,\, K_1, \lambda_1)$ satisfy Hypothesis
  \ref{hyp:condition_K}, \ref{Assumption:rhoC1} and
  \ref{hyp:U_differentiable}. Let $s=1/\alpha$ be an integer greater
  than $1$. For $F$ and $G$ cylindrical functions, for $v\in
  V_0(E_1)$, we have:
  \begin{multline*}
    \int_{\chi_\Lambda} \nabla_{v} F(\omega) G(\omega)\d\mu_{\alpha,\,
      sK_{1,\Lambda},\lambda_1}(\omega) =-\int_{\chi_\Lambda}
    F(\omega)\nabla_{v} G(\omega)\d\mu_{\alpha,\, K_{1,\,\Lambda},\lambda_1}(\omega)\\
    +\frac 1\alpha\int_{\chi_\Lambda}
    F(\omega)G(\omega)\left(\sum_{\eta\subset\omega}
      (B^{\lambda_1}_v(\eta)+\nabla_vU(\eta))R(\eta,\, \omega)
    \right)\d\mu_{\alpha,\,sK_{1,\Lambda},\lambda_1}(\omega).
  \end{multline*}
\end{Corollary}

\subsection{Consequences}
\label{sec:consequences}

We define the norm $||.||_{2,1}$ on
$\mathcal{F}C^{\infty}_b(\mathcal{D},\chi)$ by:
\begin{align*}
  ||F||^2_{2,1}&=||F||^2_{L^2(\mu)}+\E\left[||\nabla F||^2\right]\\
  &=\E\left[F^2\right]+\E\left[\int |\nabla_x F|^2\d\xi(x)\right].
\end{align*}
and we call $\mathcal{D}_{2,1}$ the closure of
$\mathcal{F}C^{\infty}_b(\mathcal{D},\chi)$ for the norm
$||.||_{2,1}$. A trivial consequence of the previous results is that,
for any $\alpha$-DPPP known to exist, the operator $\nabla$ is
closable and can thus be extended to $\mathcal{D}_{2,1}.$ With the
same lines of proof we retrieve the result of (\cite{MR2209150}),
which says that the Dirichlet form:
\begin{math}
  \mathcal{E}(F,F)=\esp{}{\langle\nabla F , \nabla F\rangle}
\end{math}
is closable.

\providecommand{\bysame}{\leavevmode\hbox to3em{\hrulefill}\thinspace}
\providecommand{\MR}[1]{}

\end{document}